\documentclass{amsart}
\usepackage{graphicx}
\usepackage{amsfonts}
\usepackage{eucal}
\usepackage{amscd}
\usepackage{amssymb}
\usepackage{xypic}
\xyoption{all}
\DeclareFontFamily{OT1}{pzc}{}
\DeclareFontShape{OT1}{pzc}{m}{it}%
             {<-> s * [1.200] pzcmi7t}{}
\DeclareMathAlphabet{\mathscr}{OT1}{pzc}%
                                 {m}{it}
\newcommand{\Mod}[1]{\mbox{ (mod }#1 \mbox{)}}

\newcommand{\Cat}[1]{ {\mathfrak{#1}} }

\newcommand{\alg}[1]{\boldsymbol{#1}}
\newcommand{\scheme}[1]{\boldsymbol{\underline{#1}}}

\newcommand{\ff}[1]{\mathbf{#1}}

\newcommand{\Sp}[1]{\mathscr{#1}}

\newcommand{\Hom}{\mathsf{Hom} }

\newcommand{\ints}{\mathbb Z}

\newcommand{\reals}{\mathbb R}
\newcommand{\complex}{\mathbb C}

\newcommand{\adeles}{\mathbb A}

\newcommand{\OO}{\mathcal{O}}

\newcommand{\Weil}{\mathcal{W}}

\newcommand{\isom}{\cong}
\newcommand{\ident}{\equiv}

\newcommand{\TwoTerm}[3]{\xymatrix{{#1} \ar[r]^{#2} & {#3}}}
\newtheorem{thm}{Theorem}[section]
\newtheorem{cor}[thm]{Corollary}
\newtheorem{lem}[thm]{Lemma}
\newtheorem{prop}[thm]{Proposition}
\theoremstyle{definition}
\newtheorem{defn}[thm]{Definition}
\newtheorem{rem}[thm]{Remark}
\numberwithin{equation}{section}

\begin{document}

\title{Metaplectic Tori over Local Fields}%
\author{Martin H. Weissman}%
\address{Dept. of Mathematics, University of California, Santa Cruz, CA 95064}%
\email{weissman@ucsc.edu}
%\thanks{}%
%\keywords{metaplectic, torus, local, Langlands}
%\date{}%
%\dedicatory{}%
%\commby{}%
% ----------------------------------------------------------------
\begin{abstract}
Smooth irreducible representations of tori over local fields have been parameterized by Langlands, using class field theory and Galois cohomology.  This paper extends this parameterization to some central extensions of such tori, which arise naturally in the setting of nonlinear covers of reductive groups.
\end{abstract}
\maketitle
% ----------------------------------------------------------------
\tableofcontents

\section{Introduction}

\subsection{Motivation}

Let $\alg{T}$ be an algebraic torus over a local field $F$; let $T = \alg{T}(F)$.  Let $L/F$ be a finite Galois extension over which $\alg{T}$ splits, with $\Gamma = Gal(L/F)$.  Let $\Sp{X}(T)$ denote the group of continuous characters of $T$ with values in $\complex^\times$.  In a preprint from 1968, now appearing as \cite{Lan} (cf. the 1985 article by Labesse \cite{Lab}), Langlands proves the following:
\begin{thm}
There is a natural isomorphism:
$$\Sp{X}(T) \isom H_c^1(\Weil_{L/F}, \Sp{\hat T}),$$
where $\Weil_{L/F}$ denotes the Weil group of $L/F$, and $\Sp{\hat T}$ denotes the complex dual torus of $\alg{T}$.
\end{thm}

We may consider $\alg{T}$ as a sheaf of groups, on the big Zariski site over $F$.  In addition, we may consider $\alg{K}_2$ as such a sheaf, using Quillen's algebraic K-theory.  Let $\alg{T}'$ be a central extension of $\alg{T}$ by $\alg{K}_2$, in the category of sheaves of groups on the big Zariski site over $F$.  Such objects are introduced and studied extensively by Deligne and Brylinski in \cite{D-B}.

Let $T' = \alg{T}'(F)$ be the resulting extension of $T$ by $K_2 = \alg{K}_2(F)$.  If $F \not \isom \complex$ and $F$ has sufficiently many $n^{th}$ roots of unity, one may push forward the central extension $T'$ via the Hilbert symbol to obtain a central extension $\tilde T$:
$$1 \rightarrow \mu_n \rightarrow \tilde T \rightarrow T \rightarrow 1.$$

We are interested in the set $\Sp{Irr}_\epsilon(\tilde T)$ of irreducible genuine representations of $\tilde T$, as defined in Section \ref{GR}.  Such representations arise frequently in the literature on ``metaplectic groups'', especially when considering principal series representations of nonlinear covers of reductive groups (cf. \cite{Sav}, \cite{K-P}, \cite{ABPTV}, among others).  It is the goal of this paper to parameterize the set $\Sp{Irr}_\epsilon(\tilde T)$ in a way which naturally generalizes the aforementioned theorem of Langlands.

\subsection{Main Results}

Associated to the central extension $\alg{T}'$, Deligne and Brylinski associate two functorial invariants:  an integer-valued quadratic form $Q$ on the cocharacter lattice $Y$ of $\alg{T}$, and a $\Gamma$-equivariant central extension $\tilde Y$ of $Y$ by $L^\times$.  Associated to $Q$, one has a symmetric bilinear form $B_Q \colon Y \otimes_\ints Y \rightarrow \ints$.

Define:
$$Y^\# = \{ y \in Y \mbox{ such that } B_Q(y,y') \in n \ints \mbox{ for all } y' \in Y \}.$$
Similarly, define:
$$Y^{\Gamma \#} = \{ y \in Y \mbox{ such that } B_Q(y,y') \in n \ints \mbox{ for all } y' \in Y^\Gamma \}.$$
Associated to the inclusion $\iota \colon Y^\# \hookrightarrow Y$, there is an isogeny of complex tori:  $\hat \iota \colon \Sp{\hat T} \rightarrow \Sp{\hat T^\#}$.  This isogeny is also a morphism of $\Weil_{L/F}$-modules.  Associated to the sequence of inclusions $Y^\# \subset Y^{\Gamma \#} \subset Y$, there are also $F$-isogenies of $F$-tori:  $\alg{T}^\# \rightarrow \alg{T}^{\Gamma \#} \rightarrow \alg{T}$.  The main results of this paper are Theorems \ref{ST1}, \ref{TT}, and \ref{RT}.  Putting these theorems together yields the following:
\begin{thm}
Suppose that one of the following conditions is satisfied:
\begin{enumerate}
\item
$\alg{T}$ is a split torus.
\item
$F$ is nonarchimedean with residue field $\ff{f}$, $\alg{T}$ splits over an unramified extension of $F$, and $n$ is relatively prime to the characteristic of $\ff{f}$.
\item
$F \isom \reals$.
\end{enumerate}
Then, there exists a finite-to-one map:
$$\Phi \colon \Sp{Irr}_\epsilon(\tilde T) \rightarrow { {H_c^1(\Weil_{L/F}, \Sp{\hat T})} \over {H_c^1(\Weil_{L/F}, \Sp{\hat T} \rightarrow \Sp{\hat T^\#})} },$$
which intertwines canonical actions of $H^1(\Weil_{L/F}, \Sp{\hat T})$.  The finite fibres of this map are torsors for a finite group $\Sp{X}(P^\dag) = Hom(P^\dag, \complex^\times)$.  In the three cases described above, the ``packet group'' $P^\dag$ can be described by the following three conditions, respectively:
\begin{enumerate}
\item
$P^\dag$ is trivial.
\item
$$P^\dag = { {Im(\alg{\bar T}^{\Gamma \#}(\ff{f}) \rightarrow \alg{\bar T}(\ff{f}))} \over {Im(\alg{\bar T}^{\#}(\ff{f}) \rightarrow \alg{\bar T}(\ff{f}))} }.$$
\item
$$P^\dag = { {Im(\pi_0 \alg{T}^{\Gamma \#}(\reals) \rightarrow \pi_0 \alg{T}(\reals))} \over {Im(\pi_0 \alg{T}^\#(\reals) \rightarrow \pi_0 \alg{T}(\reals))} }.$$
\end{enumerate}
\end{thm}

The parameterization $\Phi$ of irreducible genuine representations is not unique; rather, it depends upon the choice of a base-point.  The choice of this base-point is a significant problem.  We identify a natural class of ``pseudo-spherical'' representations (following previous authors such as \cite{Sav} and \cite{ABPTV}).  Moreover, we parameterize pseudo-spherical irreducible representations as a torsor for a complex algebraic torus in Section \ref{PS}; perhaps more naturally, the category of pseudo-spherical representations can be identified with the category of modules over a ``quantum dual torus''.

\subsection{Acknowledgements}

We would like to thank Jeffrey Adams and Gordan Savin, for providing some advice and insight related to this paper.  In addition, we are thankful for the hospitality and excellent working environment provided by the Hausdorff Institute for Mathematics in Bonn, Germany, during the preparation of this paper.  In addition, we thank the University of Michigan for their hospitality, while this paper was being finished.  We thank Brian Conrad and Stephen DeBacker for helpful conversations at Michigan during this time.

We are heavily indebted to Pierre Deligne, who read an early draft of this paper, and gave extensive helpful comments.  His generosity has led to significant improvements in presentation and content.

% ----------------------------------------------------------------

\section{Background}
\subsection{Fields and sheaves}

$F$ will always denote a local field.  $F_{Zar}$ will denote the big Zariski site over $F$.  By this, we mean that $F_{Zar}$ is the full subcategory of the category of schemes over $F$, whose objects are schemes of finite type over $F$, endowed with the Zariski topology.  $\Cat{Set}_F$ will denote the topos of sheaves of sets over $F_{Zar}$, and $\Cat{Gp}_F$ will denote the topos of sheaves of groups over $F_{Zar}$.

Any scheme or algebraic group over $F$ will be identified with its functor of points, i.e., the associated object of $\Cat{Set}_F$ or $\Cat{Gp}_F$, respectively.  Quillen's K-theory, defined in \cite{Qui}, yields sheaves $\alg{K}_n$ of abelian groups on $F_{Zar}$.  We only work with $\alg{K}_1$ and $\alg{K}_2$, viewed as objects of $\Cat{Gp}_F$.

For any field $L$, the group $\alg{K}_2(L)$ is identified as a quotient:
$$\alg{K}_2(L) = { {L^\times \otimes_\ints L^\times} \over { \langle x \otimes (1 - x) \rangle_{1 \neq x \in L^\times} } }.$$
If $l_1, l_2 \in L$, and $l_1, l_2 \not \in \{0,1 \}$, then we write $\{ l_1, l_2 \}$ for the image of $l_1 \otimes l_2$ in $\alg{K}_2(L)$.  The bilinear form $\{ \cdot, \cdot \}$ is called the {\em universal symbol}; it is skew-symmetric.  It is usually not alternating, but $\{x,x \} = \{x, -1 \}$ for all $x \in L^\times$.

\subsection{Local Nonarchimedean Fields}

Suppose that $F$ is a nonarchimedean local field.   Then $\OO_F$ will denote the valuation ring of $F$, and $\ff{f}$ the residue field of $\OO_F$.  We let $p$ denote the characteristic of $\ff{f}$, and assume that the value group of $F$ is $\ints$.  We let $q$ denote the cardinality of $\ff{f}$.

There is a canonical short exact sequence of abelian groups, given by inclusion and valuation:
$$1 \rightarrow \OO_F^\times \rightarrow F^\times \rightarrow \ints \rightarrow 1.$$
It is sometimes convenient to split this sequence of abelian groups, by choosing a uniformizing element $\varpi \in F^\times$.  However, none of our main results depend on the which uniformizing element is chosen.

Reduction yields another canonical short exact sequence:
$$1 \rightarrow \OO_F^{\times 1} \rightarrow \OO_F^\times \rightarrow \ff{f}^\times \rightarrow 1.$$
This sequence is split by the Teichmuller lifting:
$$\Theta \colon \ff{f}^\times \rightarrow \OO_F^\times.$$

\subsection{The Weil group}
We let $\Weil_F$ denote a Weil group of $F$ as discussed by Tate in \cite{Tat}. In particular, we follow Tate's choices, and normalize the reciprocity isomorphism of nonarchimedean local class field theory, $rec \colon F^\times \rightarrow \Weil_F^{ab}$ in such a way that uniformizing elements of $F^\times$ act as the {\em geometric} Frobenius via $rec$.

When $L$ is a finite Galois extension of $F$, we continue to follow Tate \cite{Tat}, and define:
$$\Weil_{L/F} = \Weil_F / \overline{ [ \Weil_L, \Weil_L ] }.$$
There is then a short exact sequence:
$$1 \rightarrow L^\times \rightarrow \Weil_{L/F} \rightarrow Gal(L/F) \rightarrow 1.$$

\subsection{The Hilbert symbol}
We say that $F$ has {\em enough} $n^{th}$ roots of unity if $\alg{\mu}_n(F)$ has $n$ elements.  When $F$ has enough $n^{th}$ roots of unity, and $F \not \isom \complex$, the Hilbert symbol provides a non-degenerate skew-symmetric bilinear map:
$$(\cdot, \cdot)_{F,n} \colon { {F^\times} \over {F^{\times n}}} \otimes_\ints { {F^\times} \over {F^{\times n}}} \rightarrow \alg{\mu}_n(F).$$
In general, the Hilbert symbol is not alternating.  The Hilbert symbol factors through $\alg{K}_2(F)$, via the universal symbol.

The definition of the Hilbert symbol relies on a choice of reciprocity isomorphism in local class field theory -- this choice has been made earlier, in sending a uniformizing element of $F^\times$ to a geometric Frobenius.

If $F$ is nonarchimedean, and $(p,n) = 1$, then we say that the Hilbert symbol $(\cdot, \cdot)_{F,n}$ is tame.  If $p$ is odd, then in the tame case, $(\varpi, \varpi)_{F,n} = (-1)^{(q-1)/n}$, for every uniformizing element $\varpi \in F^\times$.  When $p = 2$, in the tame case, $(\varpi, \varpi)_{F,n} = 1$.  When $F \isom \reals$, $(-1,-1)_{F,2} = -1$.

\subsection{Tori}
$\alg{T}$ will always denote an algebraic torus over $F$.  Let $L$ be a finite Galois extension of $F$, over which $\alg{T}$ splits, and define $\Gamma = Gal(L/F)$.  We write $X = \Hom(\alg{T}, \alg{G}_m)$ for the character group and $Y$ for the cocharacter group $\Hom(\alg{G}_m, \alg{T})$.  We view $X$ and $Y$ as finite rank free $\ints$-modules, endowed with actions of $\Gamma$.  The groups $X$ and $Y$ are in canonical $\Gamma$-invariant duality.

The dual torus $\alg{\hat T}$ is the split torus $Spec(\ints[Y])$ over $\ints$, with the resulting action of $\Gamma$.  We write $\Sp{\hat T} = \alg{\hat T}(\complex) \ident X \otimes_\ints \complex^\times$ for the resulting $\complex$-torus, also endowed with the action of $\Gamma$.

\subsection{Central Extensions of Tori by $\alg{K}_2$}
Let $\Cat{CExt}(\alg{T}, \alg{K}_2)$ be the category of central extensions of $\alg{T}$ by $\alg{K}_2$ in $\Cat{Gp}_F$.  Let $\Cat{CExt}_\Gamma(Y, L^\times)$ be the category of $\Gamma$-equivariant extensions of $Y$ by $L^\times$.

In \cite{D-B}, Deligne and Brylinski study the following category, which we call $\Cat{DB}_{\alg{T}}$.  Its objects are pairs $(Q,\tilde Y)$, where:
\begin{itemize}
\item
$Q \colon Y \rightarrow \ints$ is a $\Gamma$-invariant quadratic form.
\item
$\tilde Y$ is a $\Gamma$-equivariant central extension of $Y$ by $L^\times$.
\item
The resulting commutator map $C \colon \bigwedge^2 Y \rightarrow L^\times$ satisfies:
$$C(y_1, y_2) = (-1)^{B_Q(y_1, y_2)}, \mbox{ for all } y_1, y_2 \in Y,$$
where $B_Q$ is the symmetric bilinear form associated to $Q$.
\end{itemize}

If $(Q_1, \tilde Y_1)$, and $(Q_2, \tilde Y_2)$ are two objects of $\Cat{DB}_{\alg{T}}$, then a morphism from $(Q_1, \tilde Y_1)$ to $(Q_2, \tilde Y_2)$ exists only if $Q_1 = Q_2$, in which case the morphisms of $\Cat{DB}_{\alg{T}}$ are the just those from $\tilde Y_1$ to $\tilde Y_2$ in $\Cat{CExt}_\Gamma(Y, L^\times)$.

In Section 3.10 of \cite{D-B}, Deligne and Brylinski construct an equivalence of categories, from $\Cat{CExt}(\alg{T}, \alg{K}_2)$ to the category $\Cat{DB}_{\alg{T}}$.

In particular, given a central extension $\alg{T}'$ of $\alg{T}$ by $\alg{K}_2$, the work of \cite{D-B} (in part following Esnault \cite{Esn}) yields a quadratic form $Q \colon Y \rightarrow \ints$, and a central extension $\tilde Y$ of $Y$ by $L^\times$.  Consider the central extension $\alg{T}'(L)$ below:
$$1 \rightarrow \alg{K}_2(L) \rightarrow \alg{T}'(L) \rightarrow \alg{T}(L) \rightarrow 1.$$
In \cite{D-B}, it is shown that the resulting commutator $C_L \colon \bigwedge^2 \alg{T}(L) \rightarrow \alg{K}_2(L)$ satisfies:
$$C_L(y_1(l_1), y_2(l_2)) = \{ l_1, l_2 \}^{B_Q(y_1, y_2)},$$
for all $y_1, y_2 \in Y$, $l_1, l_2 \in L^\times$.

\subsection{Locally Compact Abelian Groups}
An {\em LCA group} is a locally compact Hausdorff separable abelian topological group.  We work here in the category $\Cat{LCAb}$ whose objects are LCA groups, and whose morphisms are continuous homomorphisms.  Suppose that we are given a short exact sequence of LCA groups and continuous homomorphisms:
$$0 \rightarrow A \rightarrow B \rightarrow C \rightarrow 0.$$
Given a fourth LCA group $D$, the functor $\Hom(\bullet, D)$ is left-exact, yielding an exact sequence:
$$0 \rightarrow \Hom(C, D) \rightarrow \Hom(B, D) \rightarrow \Hom(A,D).$$

\subsubsection{Continuous characters}
When $A$ is an LCA group, we write $\Sp{X}(A)$ for the group of continuous homomorphisms from $A$ to the LCA group $\complex^\times$, under pointwise multiplication.  We call elements of $\Sp{X}(A)$ characters (or continuous characters) of $A$.  If $\chi \in \Sp{X}(A)$, and $\vert \chi(a) \vert = 1$ for all $a \in A$, then we say that $\chi$ is a unitary character.  We write $\hat A$ for the Pontrjagin dual of $A$, i.e., the set of unitary characters of $A$, with its natural topology as an LCA group.

We say that $A$ is an elementary LCA group, if $A \isom \reals^a \times \ints^b \times (\reals / \ints)^c \times F$, for some finite group $F$, and some non-negative integers $a,b,c$.  When $A$ is elementary, $\Sp{X}(A)$ has a natural structure as a complex algebraic group.  In the case above, $\Sp{X}(A) \isom \complex^a \times (\complex^\times)^b \times \ints^c \times \hat F$.

If $A$ is generated by a compact neighborhood of the identity, then $A$ is canonically isomorphic to the inverse limit of its elementary quotients by compact subgroups.  In this case, $\Sp{X}(A)$ is endowed with the (inductive limit) structure of a complex algebraic group.  In this paper, all LCA groups will be generated by a compact neighborhood of the identity, and thus $\Sp{X}(A)$ will be viewed as a complex algebraic group.

\subsubsection{Exactness Criteria}
Given a short exact sequence
$$0 \rightarrow A \rightarrow B \rightarrow C \rightarrow 0,$$
there are two important cases, in which the induced map $\Sp{X}(B) \rightarrow \Sp{X}(A)$ is surjective, leading to an exact sequence:
$$0 \rightarrow \Sp{X}(C) \rightarrow \Sp{X}(B) \rightarrow \Sp{X}(A) \rightarrow 0.$$
\begin{prop}
\label{EX1}
Suppose that $A$ is compact.  Then $\Sp{X}(B) \rightarrow \Sp{X}(A)$ is surjective.
\end{prop}
\proof
When $A$ is compact, every continuous character of $A$ is unitary.  The exactness of Pontrjagin duality implies that every unitary character of $A$ extends to a unitary character of $B$.  Hence $\Sp{X}(B)$ surjects onto $\Sp{X}(A)$.
\qed

\begin{prop}
\label{EX2}
Suppose that the map from $A$ to $B$ is an open embedding.  Then $\Sp{X}(B) \rightarrow \Sp{X}(A)$ is surjective.
\end{prop}
\proof
The proof is not difficult -- it follows directly from Proposition 3.3 of \cite{H-S}, for example.
\qed

\subsection{Complex varieties and groups}

We use a script letter, such as $\Sp{M}$, to denote the (complex) points of a complex algebraic variety.  It is unnecessary for us to distinguish between complex varieties and their complex points.  If $R$ is a commutative reduced finitely-generated $\complex$-algebra, then we write $\Sp{M} = \Sp{MSpec}(R)$ for the maximal ideal spectrum of $A$, viewed as a complex variety.  We view $\complex^\times$ as a complex algebraic variety, identifying:
$$\complex^\times \ident \Sp{MSpec}(\complex[\ints]),$$
where $\complex[\ints]$ denotes the group ring.  We view $\complex$ itself as an algebraic variety (the affine line over the field $\complex$).

Suppose that $\Sp{G}$ is a complex algebraic group, in other words, a group in the category of complex algebraic varieties.  A $\Sp{G}$-variety is a complex algebraic variety $\Sp{M}$, endowed with an action $\Sp{G} \times \Sp{M} \rightarrow \Sp{M}$ which is complex-algebraic.  A $\Sp{G}$-torsor is a $\Sp{G}$-variety $\Sp{M}$, such that the induced map $\Sp{G} \times \Sp{M} \rightarrow \Sp{M} \times \Sp{M}$ sending $(g,m)$ to $(g \cdot m, m)$, is an isomorphism of complex algebraic varieties.

If $\Sp{M}_1$ and $\Sp{M}_2$ are $\Sp{G}$-varieties, then a morphism of $\Sp{G}$-varieties is a complex algebraic map from $\Sp{M}_1$ to $\Sp{M}_2$, which intertwines the action of $\Sp{G}$.  Morphisms of torsors are defined in the same way.

\section{Genuine Representations of Metaplectic Tori}
\label{GR}
In this section, the following will be fixed:
\begin{itemize}
\item
$F$ will be a local field, with $F \not \isom \complex$.  $n$ will be a positive integer, such that $F$ has enough $n^{th}$ roots of unity.
\item
$\alg{T}$ will be a torus over a local field $F$, which splits over a finite Galois extension $L/F$, with $\Gamma = Gal(L/F)$.  $X$ and $Y$ will be the resulting character and cocharacter groups.
\item
$\alg{T}'$ will be an extension of $\alg{T}$ by $\alg{K}_2$ in $\Cat{Gp}_F$.
\item
$(Q,\tilde Y)$ will be the Deligne-Brylisnki invariants of $\alg{T}'$.  $B$ will be the symmetric bilinear form associated to $Q$.
\item
$\epsilon \colon \alg{\mu}_n(F) \rightarrow \complex^\times$ will be a fixed injective character.
\end{itemize}

\subsection{Heisenberg Groups}

Suppose that $S$ is an LCA group, and $A$ is a finite cyclic abelian group endowed with a faithful unitary character $\epsilon \colon A \rightarrow \complex^\times$. Suppose that $\tilde S$ is a locally compact group, which is a central extension of $S$ by $A$ (in the category of locally compact groups and continuous homomorphisms:
$$1 \rightarrow A \rightarrow \tilde S \rightarrow S \rightarrow 1.$$
In this situation, the commutator on $\tilde S$ descends to a unique alternating form:
$$C \colon \bigwedge^2 S \rightarrow A.$$
Let $Z(\tilde S)$ be the center of $\tilde S$.  Then $Z(\tilde S)$ is the preimage of a subgroup $Z^\dag(S) \subset S$, where:
$$Z^\dag(S) = \{ s \in S \mbox{ such that } C(s,s') = 1 \mbox{ for all } s' \in S \}.$$
Throughout this paper, the following condition will be satisfied, and hence we assume that:
\begin{quote}
$Z^\dag(S)$ is an open subgroup of finite index in $S$.
\end{quote}

We define the following two sets:
\begin{itemize}
\item
The set $\Sp{X}_\epsilon(\tilde S)$ of continuous genuine characters of $\tilde S$.  These are elements of $\Sp{X}(\tilde S)$ whose restriction to $A$ equals $\epsilon$.
\item
The set $\Sp{Irr}_\epsilon(\tilde S)$ of irreducible genuine representations of $\tilde S$.  These are irreducible (algebraic) representations $(\pi, V)$ of $\tilde S$ on a complex vector space, on which $Z(\tilde S)$ acts via a continuous genuine character.  In particular, since $Z(\tilde S)$ will always have finite index in $\tilde S$, these are finite-dimensional representations.
\end{itemize}

We often use the following analogue of the Stone von-Neumann theorem:
\begin{thm}
\label{SvN}
Suppose that $\chi \in \Sp{X}_\epsilon(Z(\tilde S))$ is a genuine continuous character.  Let $\tilde M$ denote a maximal commutative subgroup of $\tilde S$.  Then there exists an extension $\tilde \chi \in \Sp{X}(\tilde M)$ of $\chi$ to $\tilde M$.  Define a representation of $\tilde S$ by:
$$(\pi_\chi, V_\chi) = Ind_{\tilde M}^{\tilde S} \tilde \chi.$$
Algebraic induction suffices here, since we always assume that $Z(\tilde S)$ has finite index in $\tilde S$.  Then we have the following:
\begin{enumerate}
\item
The representation $(\pi_\chi, V_\chi)$ is irreducible.
\item
The representation $(\pi_\chi, V_\chi)$ has central character $\chi$.
\item
The isomorphism class of $(\pi_\chi, V_\chi)$ depends only upon $\chi$, and not upon the choices of subgroup $\tilde M$ and extension $\tilde \chi$.
\item
Every irreducible representation of $\tilde S$, on which $Z(\tilde S)$ acts via $\chi$, is isomorphic to $(\pi_\chi, V_\chi)$.
\end{enumerate}
\end{thm}
\proof
Extension of $\chi$ to $\tilde M$ follows from Proposition \ref{EX2}.  All but the last claim are proven in Section 0.3 of \cite{K-P}, and follow directly from Mackey theory.  The last claim follows from the previous claims and Frobenius reciprocity.
\qed

\subsection{Metaplectic Tori over Local Fields}

The central extension of $\alg{T}$ by $\alg{K}_2$ yields a central extension of groups:
$$1 \rightarrow \alg{K}_2(F) \rightarrow \alg{T}'(F) \rightarrow \alg{T}(F) \rightarrow 1.$$
Since $F$ is assumed to have enough $n^{th}$ roots of unity, the Hilbert symbol allows us to push forward this extension to get:
$$1 \rightarrow \mu_n \rightarrow \tilde T \rightarrow T \rightarrow 1,$$
where $\mu_n = \alg{\mu}_n(F)$, and $T = \alg{T}(F)$.  By results of Sections 10.2, 10.3 of \cite{D-B}, following Moore \cite{Moo}, this is a topological central extension, of the LCA group $T$ by the LCA group $\mu_n$.

In this case, the center $Z(\tilde T)$ has finite index in $\tilde T$.  Furthermore, Theorem \ref{SvN} implies that:
\begin{prop}
There is a natural bijection between the set $\Sp{Irr}_\epsilon(\tilde T)$ of irreducible genuine representations of $\tilde T$, and the set $\Sp{X}_\epsilon(Z(\tilde T))$ of genuine characters of $Z(\tilde T)$.
\end{prop}

There is a short exact sequence of LCA groups:
$$1 \rightarrow \mu_n \rightarrow Z(\tilde T) \rightarrow Z^\dag(T) \rightarrow 1.$$
From Proposition \ref{EX1}, it follows that:
\begin{prop}
The space $\Sp{X}_\epsilon(Z(\tilde T))$ of genuine continuous characters of $Z(\tilde T)$ is a $\Sp{X}(Z^\dag(T))$-torsor.
\end{prop}
As an immediate consequence, we find:
\begin{cor}
The set $\Sp{Irr}_\epsilon(\tilde T)$ is a $\Sp{X}(Z^\dag(T))$-torsor.
\end{cor}
In particular, we give $\Sp{Irr}_\epsilon(\tilde T)$ the structure of a complex algebraic variety, in such a way that it is a complex algebraic $\Sp{X}(Z^\dag(T))$-torsor.

Since $Z^\dag(T)$ is a finite index subgroup of $T$, restriction of continuous characters yields a surjective homomorphism of complex algebraic groups:
$$res \colon \Sp{X}(T) \rightarrow \Sp{X}(Z^\dag(T)).$$
As a result, the set $\Sp{Irr}_\epsilon(\tilde T)$ is a homogeneous space for $\Sp{X}(T)$, or equivalently (by Langlands theorem \cite{Lan}), a homogeneous space for $H_c^1(\Weil_{L/F}, \Sp{\hat T})$.

\section{Split Tori}
In this section, we carry on the assumptions of the previous section.  In addition, we assume that $\alg{T}$ is a split torus of rank $r$ over $F$.  Thus, there is a canonical identification $\alg{T}(F) \ident Y \otimes_\ints F^\times$.  We are interested in parameterizing $\Sp{Irr}_\epsilon(\tilde T)$.  By the results of the previous section, we may describe this set, up to a choice of base point, by describing the set $Z^\dag(T)$.

\subsection{An isogeny}

Recall that $B \colon Y \otimes_\ints Y \rightarrow \ints$ is the symmetric bilinear form associated to $Q$.  It allows us to construct a subgroup of finite index $Y^\# \subset Y$:
$$Y^\# = \{ y \in Y \mbox{ such that } B(y,y') \in n \ints \mbox{ for all } y' \in Y \}.$$
Note that we suppress mention of $Q$, $B$, and $n$, in our notation $Y^\#$.

The subgroup $Y^\#$ can be related to the ``Smith normal form'' of the bilinear form $B$.  Namely, there exists a pair of group isomorphisms $\alpha, \beta$:
$$\xymatrix{ \ints^r & Y  \ar[l]_\beta \ar[r]^\alpha & \ints^r },$$
such that for all $y_1, y_2 \in Y$, one has:
$$B_Q(y_1, y_2) = D(\alpha(y_1), \beta(y_2)),$$
and $D$ is a symmetric bilinear form on $\ints^r$ represented by a diagonal matrix with entries $(d_1, \ldots, d_r)$ (the elementary divisors).  Let $e_j$ denote the smallest positive integer such that $d_j e_j \in n \ints$, for every $1 \leq j \leq r$.  Then we find that:
$$Y^\# = \alpha^{-1}(e_1 \ints \oplus e_2 \ints \oplus \cdots \oplus e_r \ints).$$

Let $\iota \colon Y^\# \rightarrow Y$ denote the inclusion of $\ints$-modules.  Since $Y^\#$ has finite index in $Y$, this corresponds to an $F$-isogeny of split tori:
$$\alg{\iota} \colon \alg{T}^\# \rightarrow \alg{T},$$
where $\alg{T}^\#$ is the split algebraic torus with cocharacter lattice $Y^\#$.  From the previous observations, we find that:
$$\iota(T^\#) = \iota(\alg{T}^\#(F)) = \alpha^{-1} \left( F^{\times e_1} \times \cdots \times F^{\times e_r} \right).$$

\subsection{Describing the Center}
Recall that extension $\tilde T$ of $T$ by $\mu_n$ yields a commutator:
$$C \colon \bigwedge^2 T \rightarrow \mu_n.$$
This commutator can be directly related (by \cite{D-B}) to the bilinear form $B$.  If $u_1,u_2 \in F^\times$, and $y_1, y_2 \in Y$, then one may directly compute:
$$C(y_1(u_1), y_2(u_2)) = (u_1, u_2)_n^{B(y_1, y_2)}.$$
The diagonalization of $B$, via group isomorphisms $\alpha, \beta$, yields two isomorphisms of $F$-tori:
$$\xymatrix{ \alg{G}_m^r  & \alg{T} \ar[l]_{\alg{\beta}} \ar[r]^{\alg{\alpha}} & \alg{G}_m^r },$$

One arrives at a bilinear form on $(F^\times)^r$, given by:
$$\Delta(\vec z_1, \vec z_2) = \prod_{j=1}^r \left( z_1^{(j)}, z_2^{(j)} \right)_n^{d_j}.$$
This is related to the commutator $C$ by:
$$C(t_1, t_2) = \Delta(\alpha(t_1), \beta(t_2)).$$

We can now characterize $Z^\dag(T)$:
\begin{prop}
The subgroup $Z^\dag(T)$ of $T$ is equal to the image of the isogeny $\iota$ on the $F$-rational points:
$$Z^\dag(T) = \iota(T^\#).$$
\label{CSMT}
\end{prop}
\proof
We find that:
\begin{eqnarray*}
t_1 \in Z^\dag(T) & \Leftrightarrow & C(t_1, t_2) = 1  \mbox{ for all } t_2 \in T, \\
& \Leftrightarrow & \Delta(\alpha(t_1), \beta(t_2)) = 1 \mbox{ for all } t_2 \in T, \\
& \Leftrightarrow & \Delta(\alpha(t_1), \vec z_2) = 1 \mbox{ for all } \vec z_2 \in (F^\times)^r, \\
& \Leftrightarrow & \alpha(t_1) \in \left( F^{\times e_1}  \times \cdots \times F^{\times e_r} \right) \\
& \Leftrightarrow & t_1 \in \iota(T^\#).
\end{eqnarray*}
The penultimate step above follows from the non-degeneracy of the Hilbert symbol.
\qed

We have now proven:
\begin{thm}
If $F$ is a local field, then $\Sp{Irr}_\epsilon(\tilde T)$ is a torsor for $\Sp{X}(\iota(T^\#))$.
\end{thm}

\subsection{Character groups}

The previous theorem motivates the further analysis of the group $\Sp{X}(\iota(T^\#))$.  We write $\iota^\ast$ for the pullback homomorphism:
$$\iota^\ast \colon \Sp{X}(T) \rightarrow \Sp{X}(T^\#).$$
\begin{prop}
There is a natural identification:
$$\Sp{X}(\iota(T^\#)) \ident Im(\iota^\ast).$$
\end{prop}
\proof
There are short exact sequences of LCA groups:
$$1 \rightarrow Ker(\iota) \rightarrow T^\# \rightarrow \iota(T^\#) \rightarrow 1,$$
$$1 \rightarrow \iota(T^\#) \rightarrow T \rightarrow Cok(\iota) \rightarrow 1,$$
Using Propositions \ref{EX1} and \ref{EX2}, we arrive at short exact sequences of character groups:
$$1 \rightarrow \Sp{X}(\iota(T^\#) \rightarrow \Sp{X}(T^\#) \rightarrow \Sp{X}(Ker(\iota)) \rightarrow 1,$$
$$1 \rightarrow \Sp{X}(Cok(\iota)) \rightarrow \Sp{X}(T) \rightarrow \Sp{X}(\iota(T^\#)) \rightarrow 1,$$
Since $\Sp{X}(T)$ surjects onto $\Sp{X}(\iota(T^\#))$, we find that the image of $\iota^\ast \colon \Sp{X}(T) \rightarrow \Sp{X}(T^\#)$ is equal to the image of the injective map $\Sp{X}(\iota(T^\#)) \rightarrow \Sp{X}(T^\#)$.
\qed

\subsection{The Dual Complex}

The isogeny of split $F$-tori $\alg{\iota} \colon \alg{T}^\# \rightarrow \alg{T}$ yields an isogeny of the complex dual tori:
$$\hat \iota \colon \Sp{\hat T} \rightarrow \Sp{\hat T}^\#.$$

One may pull back continuous characters:
$$\iota^\ast \colon \Sp{X}(T) \rightarrow \Sp{X}(T^\#).$$

The following result follows from local class field theory, and demonstrates the naturality of Langlands' classification \cite{Lan}:
\begin{prop}
There is a commutative diagram of complex algebraic groups, whose rows are the reciprocity isomorphisms of local class field theory:
$$\xymatrix{
\Sp{X}(T) \ar[r] \ar[d]^{\iota^\ast} & H_c^1(\Weil_F, \Sp{\hat T}) \ar[d]^{\hat \iota} \\
\Sp{X}(T^\#) \ar[r] & H_c^1(\Weil_F, \Sp{\hat T}^\#)
}$$
\end{prop}
Note that since $\alg{T}$ and $\alg{T}^\#$ are split tori, the continuous cohomology groups are quite simple:
$$H_c^1(\Weil_F, \Sp{\hat T}) = \Hom_c(\Weil_F, \Sp{\hat T}).$$

From this result, we now find:
\begin{cor}
There is a natural identification:
$$\Sp{X}(Z^\dag(T)) \ident Im \left(\hat \iota \colon H^1(\Weil_F, \Sp{\hat T}) \rightarrow H^1(\Weil_F, \Sp{\hat T}^\#) \right).$$
\end{cor}

\subsection{Parameterization by Hypercohomology}
We may parameterize representations, using the hypercohomology of the complex of tori:
$$\TwoTerm{\Sp{\hat T}}{\hat \iota}{\Sp{\hat T}^\#}$$
concentrated in degrees zero and one.  We follow the treatment of Kottwitz-Shelstad in the Appendices of \cite{K-S}, when discussing hypercohomology of Weil groups, with coefficients in complexes of tori.  In particular, we concentrate the complexes in degrees 0 and 1 following \cite{K-S}, and not in degrees -1 and 0 as in \cite{Bor}.

There is a long exact sequence in cohomology, which includes:
$$\xymatrix{
H_c^1(\Weil_F, \Sp{\hat T} \ar[r] &  \Sp{\hat T}^\#) \ar[r]^-{\eta} & H_c^1(\Weil_F, \Sp{\hat T}) \ar[r] & H_c^1(\Weil_F, \Sp{\hat T}^\#).}$$

\begin{lem}
The homomorphism $\eta$ is injective.
\end{lem}
\proof
Extending the long exact sequence above, it suffices to prove the surjectivity of the preceding homomorphism:
$$H_c^0(\Weil_F, \Sp{\hat T}) \rightarrow H_c^0(\Weil_F, \Sp{\hat T}^\#).$$
But $\Sp{\hat T}$ and $\Sp{\hat T}^\#$ are complex tori, trivial as $\Weil_F$-modules, and $\hat \iota$ is an isogeny.  Therefore, the above map is surjective.
\qed

From this lemma, we identify the hypercohomology group $H_c^1(\Weil_F, \Sp{\hat T} \rightarrow \Sp{\hat T}^\#)$ with a subgroup of $H_c^1(\Weil_F, \Sp{\hat T})$.

\begin{lem}
The group $H_c^1(\Weil_F, \Sp{\hat T} \rightarrow \Sp{\hat T}^\#)$ is finite.
\end{lem}
\proof
Since $\hat \iota$ is an isogeny, it has finite kernel and cokernel.  There is a long exact sequence, which includes the following terms:
$$H_c^1(\Weil_F, Ker(\hat \iota)) \rightarrow H_c^1(\Weil_F, \Sp{\hat T} \rightarrow \Sp{\hat T}^\#) \rightarrow H_c^1(\Weil_F, Cok(\hat \iota)).$$
The lemma follows.
\qed

This leads to the first main theorem:
\begin{thm}
\label{ST1}
There exists an isomorphism, in the category of varieties over $\complex$ endowed with an action of $H_c^1(\Weil_F, \Sp{\hat T})$:
$$\Sp{Irr}_\epsilon(\tilde T) \isom { {H_c^1(\Weil_F, \Sp{\hat T})} \over { H_c^1(\Weil_F, \Sp{\hat T} \rightarrow \Sp{\hat T}^\#)} }.$$
\end{thm}

\begin{rem}
The global analogue of this result also seems to hold.  Letting $\Sp{Irr}_\epsilon^{aut}(\tilde T_\adeles)$ denote the appropriate set of genuine automorphic representations of $\tilde T_\adeles$, it seems likely that:
$$\Sp{Irr}_\epsilon^{aut}(\tilde T_\adeles) \isom { {H_c^1(\Weil_F, \Sp{\hat T})} \over {H_c^1(\Weil_F, \Sp{\hat T} \rightarrow \Sp{\hat T}^\#) } },$$
when $\alg{T}$ is a split torus over a global field $F$.  The proof follows the same techniques (together with the Hasse principle for isogenies of split tori), but requires some analytic care with extension of continuous characters and the appropriate Stone von-Neumann theorem.  We hope to treat this global theorem in a future paper.
\end{rem}

\section{Unramified Tori}

For split tori, the cohomology groups that arise in Theorem \ref{ST1} are quite simple, since $\Weil_F$ acts trivially on $\Sp{\hat T}$ and $\Sp{\hat T}^\#$.  In fact, the statement of this theorem makes sense, even when $\alg{T}$ is a nonsplit torus.  However, for general nonsplit tori, it seems that our explicit methods are insufficient to prove such a result.  For ``tame covers'' of unramified tori over local nonarchimedean fields, such a paramaterization is possible.

In this section, we fix the following:
\begin{itemize}
\item
$\alg{T}$ will be a a local nonarchimedean field $F$, which splits over a finite {\em unramified} Galois extension $L/F$, with $\Gamma = Gal(L/F)$.  $X$ and $Y$ will be the resulting character and cocharacter groups.
\item
We define $d = [L:F]$, and write $\ff{l}$ for the residue field of $\OO_L$ (note that $\ff{l}$ has cardinality $q^d$).
\item
We fix a uniformizer $\varpi$ of $F^\times$ (and hence of $L^\times$ as well).
\item
Let $\gamma$ be the generator of $\Gamma$, which acts upon $\ff{l}$ via $\gamma(x) = x^q$.  Let $r = (q^d - 1) / (q-1) = \# (\ff{l}^\times / \ff{f}^\times )$.
\item
$\alg{T}'$ will be an extension of $\alg{T}$ by $\alg{K}_2$ in $\Cat{Gp}_F$.
\item
$(Q,\tilde Y)$ will be the Deligne-Brylisnki invariants of $\alg{T}'$.  $B$ will be the symmetric bilinear form associated to $Q$.
\item
$n$ will be a positive integer, such that $F$ has enough $n^{th}$ roots of unity.  We also assume that $(p,n) = 1$.
\item
$\epsilon \colon \alg{\mu}_n(F) \rightarrow \complex^\times$ will be a fixed injective character.
\item
If $W$ is a subgroup of $Y$, then we write $W^\Gamma$ for the subgroup of $\Gamma$-fixed elements of $W$.  We also define:
$$W^\# = \{ y \in Y \mbox{ such that } B(y, w) \in n \ints \mbox{ for all } w \in W \}.$$
\end{itemize}

\subsection{$\ints[\Gamma]$-modules}

$\Gamma$ is a cyclic group, generated by $\gamma$, of order $d$.  Let $\ints[\Gamma]$ denote the integral group ring of $\Gamma$.  We define the following elements of $\ints[\Gamma]$:
\begin{itemize}
\item
Let $Tr = \sum_{i=0}^{d-1} \gamma^i$, and let $Tr_q = \sum_{i=0}^{d-1} q^i \gamma^i$.
\item
Let $\delta = \gamma - 1$, and let $\delta_q = q \gamma - 1$.
\end{itemize}
Note that:
$$Tr \circ \delta = 0, \mbox{ and } Tr_q \circ \delta_q = q^d - 1.$$
When $M$ is an $\ints[\Gamma]$-module, we let $\bar M = M / (q^d-1) M$.  We write $M^\Gamma$ for the $\Gamma$-invariant $\ints$-submodule of $M$.  Therefore,
$$M^\Gamma = \{ m \in M \mbox{ such that } \delta m = 0 \}.$$
We define $\bar M^{\Gamma, q}$ by:
$$\bar M^{\Gamma, q} = \{ \bar m \in \bar M \mbox{ such that } \delta_q \bar m = 0 \}.$$
\begin{prop}
Suppose that $M$ is an $\ints[\Gamma]$-module.  Then,
$$Tr(M) \subset M^\Gamma \mbox{ and } Tr_q(\bar M) \subset \bar M^{\Gamma, q}.$$
\end{prop}
\proof
The first inclusion is obvious.  For the second inclusion, suppose that $\bar m \in \bar M$.  Then,
$$\delta_q Tr_q \bar m = Tr_q \delta_q \bar m = (q^d - 1) \bar m = 0.$$
\qed
\begin{prop}
Suppose that $M$ is a $\ints[\Gamma]$-module, which is free as an $\ints$-module.  Then, $\delta_q$ and $Tr_q$ act as injective  endomorphisms of $M$, and
$$Im(\delta_q) = \{ m \in M \mbox{ such that } Tr_q m \in (q^d - 1) M \}.$$
\end{prop}
\proof
Since $M$ is free as an $\ints$-module, $\delta_q \circ Tr_q = Tr_q \circ \delta_q = q^d-1$ acts as an injective endomorphism of $M$.  Hence $\delta_q$, and $Tr_q$ must also act as injective endomorphisms of $M$, proving the first assertion.

Since $Tr_q \circ \delta_q = q^d - 1$, it follows that:
$$Im(\delta_q) \subset \{ m \in M \mbox{ such that } Tr_q m \in (q^d-1) M \}.$$
In the other direction, if $Tr_q m \in (q^d-1) M$, then $Tr_q m = Tr_q \delta_q m'$, for some $m' \in M$.  Since $Tr_q$ acts via an injective endomorphism, it follows that $m = \delta_q m'$.
\qed

\subsection{Unramified Tori}

Much of our treatment of {\em unramified tori} is inspired by Section 2 of Ono \cite{Ono}.  Recall that $X$ and $Y$ are naturally $\ints[\Gamma]$-modules, in such a way that the pairing is $\Gamma$-invariant.

We fix a smooth model $\scheme{T}$ of $\alg{T}$ over $\OO_F$.  We make the following identifications:
$$T_L = \alg{T}(L) \ident Y \otimes_\ints L^\times, \mbox{ and } T_F = \alg{T}(F) \ident (Y \otimes_\ints L^\times)^\Gamma.$$
Similarly, for the integral points, we identify:
$$T_L^\circ = \scheme{T}(\OO_L) \ident Y \otimes_\ints \OO_L^\times, \mbox{ and } T_F^\circ = \scheme{T}(\OO_F) \ident (Y \otimes_\ints \OO_L^\times)^\Gamma.$$
We write $\alg{\bar T}$ for the special fibre of $\scheme{T}$.  Then, we also identify:
$$\bar T_{\ff{l}} = \alg{\bar T}(\ff{l})\ident Y \otimes_\ints \ff{l}^\times, \mbox{ and } \bar T_{\ff{f}} = \alg{\bar T}(\ff{f}) \ident (Y \otimes_\ints \ff{l}^\times)^\Gamma.$$
There are natural reduction homomorphisms:
$$T_L^\circ \rightarrow \bar T_{\ff{l}}, \mbox{ and } T_F^\circ \rightarrow \bar T_{\ff{f}}.$$
Let $T_L^1$ and $T_F^1$ denote the kernels of these reduction maps.  The reduction morphisms are split by the Teichmuller lift, and we arrive at a decomposition of $\ints[\Gamma]$-modules:
$$T_L^\circ \ident T_L^1 \times \bar T_{\ff{l}}.$$
Together with the valuation map, we arrive at a short exact sequence of $\ints[\Gamma]$-modules:
$$1 \rightarrow T_L^1 \times \bar T_{\ff{l}} \rightarrow T_L \rightarrow Y \rightarrow 1.$$
The choice of ($\Gamma$-invariant) uniformizing element $\varpi$ splits this exact sequence, leading to a decomposition of $\ints[\Gamma]$-modules:
$$T_L \ident Y \times \bar T_{\ff{l}} \times T_L^1.$$

We use this decomposition to ``get our hands on'' elements of $T_L$.  First, every element of $T_L$ can be expressed as $y(\varpi) t^\circ$, for uniquely determined $y \in Y$, $t^\circ \in T_L^\circ$.  Let $\theta_\ff{l}$ denote a generator of the cyclic group $\ff{l}^\times$, and $\theta_\ff{f} = \theta_\ff{l}^r$.  Thus $\theta_\ff{f}$ is a generator of the cyclic group $\ff{f}^\times$.  Let $\vartheta_L \in \OO_L^\times $ and $\vartheta_F \in \OO_F^\times$ denote the Teichmuller lifts of $\theta_\ff{l}$ and $\theta_\ff{f}$, respectively.

Let $\zeta_L = (\varpi, \vartheta_L)_{L,q^d-1}$.  Let $\zeta_F = \zeta_L^r$.  Note that $\zeta_L$ is a primitive $(q^d-1)^{th}$ root of unity, and $\zeta_F$ is a primitive $(q-1)^{th}$ root of unity.

Recall that $\bar Y = Y / (q^d-1) Y$; thus, for $\bar y \in \bar Y$, it makes sense to write $\bar y(\vartheta_L)$ as an element of $T_L^\circ$.  According to the decomposition $T_L \ident Y \times \bar T_{\ff{l}} \times T_L^1$, every element $t \in T_L$ has a unique expression:
$$t = y_1(\varpi) \bar y_2(\vartheta_L) t^1,$$
where $y_1 \in Y$, $\bar y_2 \in \bar Y$, and $t^1 \in T_L^1$.
To determine when such an expression lies in $T_F$, we have the following characterization:
\begin{prop}
An element $y_1(\varpi) \bar y_2(\vartheta_L) t^1$ of $T_L$, with $y_1, \bar y_2, t^1$ as above, lies in $T_F$ if and only if the following three conditions hold:
\begin{itemize}
\item
$y_1 \in Y^\Gamma$.  In other words, $\delta(y_1) = 0$.
\item
$\bar y_2 \in \bar Y^{\Gamma, q}$.  In other words, $\delta_q(\bar y_2) = 0$.
\item
$t^1 \in T_F^1$.
\end{itemize}
\end{prop}
\proof
By the $\Gamma$-invariance of the decomposition $T_L \ident Y \times \bar T_\ff{l} \times T_L^1$, we find that $y_1(\varpi) \bar y_2(\vartheta_L) t^1 \in T_F$ if and only if the three factors are fixed by $\Gamma$.  The proposition follows from three observations:
\begin{itemize}
\item
Since $\varpi \in F$, we have $y_1(\varpi) \in T_L^\Gamma$ if and only if $y_1 \in Y^\Gamma$.
\item
Since $\gamma(\vartheta_L) = \vartheta_L^q$, we find that $\bar y_2(\vartheta_L) \in T_L^\Gamma$ if and only if:
$$\bar y_2 = q \gamma(\bar y_2), \mbox{ in } \bar Y.$$
\item
Since the reduction map intertwines the action of $\Gamma$, we have $t^1 \in (T_L^1)^\Gamma$ if and only if $t^1 \in T_F^1$.
\end{itemize}
\qed

\subsection{Tame Metaplectic Unramified Tori}

The structure of $\alg{T}'(L)$ and $\alg{T}'(F)$ is based on Sections 12.8-12.12 of \cite{D-B}.  In particular, letting $T_L' = \alg{T}'(L)$, and $T_F' = \alg{T}'(F)$, there is a natural commutative diagram:
$$\xymatrix{
1 \ar[r] & \alg{K}_2(F) \ar[r] \ar[d] & T_F' \ar[r] \ar[d] & T_F \ar[r] \ar[d] & 1 \\
1 \ar[r] & \alg{K}_2(L) \ar[r] & T_L' \ar[r] & T_L \ar[r] & 1.}$$
There is a natural action of $\Gamma$ on the bottom row, in such a way that $\alg{K}_2(F)$ maps to $\alg{K}_2(L)^\Gamma$, $T_F = T_L^\Gamma$, and $T_F'$ maps to $(T_L')^\Gamma$.  The tame symbols yield a commutative diagram, where the downward arrows arise from the functoriality of $\alg{K}_2$ and $\alg{K}_1$:
$$\xymatrix{
\alg{K}_2(F) \ar[r]^t \ar[d] & \ff{f}^\times \ar[d] \\
\alg{K}_2(L) \ar[r]^t & \ff{l}^\times }$$
The bottom row is a morphism of $\ints[\Gamma]$-modules.  Pushing forward $T_F'$ and $T_L'$ via the tame symbols yields a commutative diagram of locally compact groups, with exact rows:
$$\xymatrix{
1 \ar[r] & \ff{f}^\times \ar[r] \ar[d] & \tilde T_F^{t} \ar[r] \ar[d] & T_F \ar[r] \ar[d] & 1 \\
1 \ar[r] & \ff{l}^\times \ar[r] & \tilde T_L^{t} \ar[r] & T_L \ar[r] & 1.}$$
The downward arrows arise from the inclusion of $F$ in $L$, and of $\ff{f}$ in $\ff{l}$.  In Section 12.8 of \cite{D-B}, Deligne and Brylinski note the following:
\begin{prop}
In the commutative diagram above, the groups in the top row are precisely the $\Gamma$-invariant subgroups of the bottom row.  In other words, $\ff{f}^\times = (\ff{l}^\times)^\Gamma$, $T_F = T_L^\Gamma$, and $\tilde T_F^{t} = (\tilde T_L^{t} )^\Gamma$.
\end{prop}

We may push forward the covers further to obtain all {\em tame covers}.  Recall that $(p,n) = 1$, and $F$ has enough $n^{th}$ roots of unity.  Then, we find that $n \vert (q-1)$, and there is a natural surjective map:
$$\psi_F \colon \ff{f}^\times \rightarrow \alg{\mu}_n(F),$$
obtained by first applying the Teichmuller map (from $\ff{f}^\times$ to $\alg{\mu}_{q-1}(F)$), and then raising to the $m = (q-1)/n$ power.  Recall that $r = (q^d - 1)/(q-1)$.  One obtains a similar map:
$$\psi_L \colon \ff{l}^\times \rightarrow \alg{\mu}_{nr}(L),$$
obtained by applying the Teichmuller map (from $\ff{l}^\times$ to $\alg{\mu}_{q^d-1}(L)$) and then raising to the $m = (q-1)/n$ power.  The compatibility of these maps yields a new commutative diagram with exact rows:
$$\xymatrix{
1 \ar[r] & \alg{\mu}_n(F) \ar[r] \ar[d] & \tilde T_F \ar[r] \ar[d] & T_F \ar[r] \ar[d] & 1 \\
1 \ar[r] & \alg{\mu}_{nr}(L) \ar[r] & \tilde T_L \ar[r] & T_L \ar[r] & 1.}$$
With this construction, we say that $\tilde T_F$ is a {\em tame} metaplectic cover of $T_F$, and $\tilde T_L$ is a tame metaplectic cover of $T_L$ as well.  $\tilde T_F$ is identified as a subgroup of $\tilde T_L$.

Note that the commutator map for $\tilde T_L$ satisfies:
$$C_{L}(y_1(u), y_2(v)) = (u,v)_{L,nr}^{B(y_1, y_2)} = (u,v)_{L, q^d-1}^{m B(y_1, y_2)},$$
where $(\cdot, \cdot)_{L, nr}$ and $(\cdot, \cdot)_{L, q^d-1}$ denotes the appropriate Hilbert symbols (in this case, norm residue symbols) on $L^\times$.  The commutator on $T_F$ is simply the restriction of $C_L$; as a result, $Z^\dag(T_F) \supset Z^\dag(T_L) \cap T_F$, where the pre-image of $Z^\dag(T_F)$ is the center of $\tilde T_F$, and the pre-image of $Z^\dag(T_L)$ is the center of $\tilde T_L$.

\subsection{Computation of the center}

Recall that the set $\Sp{Irr}_\epsilon(\tilde T_F)$ is a torsor for $\Sp{X}(Z^\dag(T_F))$.  Therefore, we wish to study the group $Z^\dag(T_F)$ in more detail.  To this end, we first observe:
\begin{prop}
The group $T_L^1$ is contained in $Z^\dag(T_L)$.  Similarly, $T_F^1$ is contained in $Z^\dag(T_F)$.
\end{prop}
\proof
Since $(q^d - 1, p) = 1$, the Hilbert symbol (in this case, a norm-residue symbol) is trivial, when one of its ``inputs'' is contained in $\OO_L^1$.  Hence the commutator $C_L(\cdot, \cdot)$ is trivial when one of its inputs is contained in $T_L^1$.  Hence $T_L^1 \subset Z^\dag(T_L)$.  Since $Z^\dag(T_F) \supset Z^\dag(T_L) \cap T_F$, we find that $T_F^1 \subset Z^\dag(T_F)$ as well.
\qed

Since $T_F^1$ is contained in $Z^\dag(T_F)$, $Z^\dag(T_F)$ corresponds to a subgroup of $T_F / T_F^1$.  Our choice of uniformizing element, together with the previously mentioned splittings, yields a decomposition of $\ints[\Gamma]$-modules:
$$T_L / T_L^1 \ident Y \times \bar T_\ff{l}.$$
Namely, every element $t$ of $T_L / T_L^1$ can be represented by $y_1(\varpi) \bar y_2(\vartheta_L)$, for uniquely determined $y_1 \in Y$ and $\bar y_2 \in \bar Y$.

In order to describe $Z^\dag(T_F)$, we work with a number of subgroups of $Y$.  Recall that $Y^{\Gamma \#}$ is given by:
$$Y^{\Gamma \#} = \{ y \in Y \mbox{ such that } B(y, y') \in n \ints \mbox{ for all } y' \in Y^\Gamma \}.$$
Note that $Y^{\Gamma \#} \supset Y^\#$.  Also, it is important to distinguish between $Y^{\Gamma \#}$ and $Y^{\# \Gamma} = (Y^\#)^\Gamma$.

\begin{lem}
There are inclusions of $\ints[\Gamma]$-modules, of finite index in $Y$:
$$Y \supset Y^{\Gamma \#} \supset Y^\# \supset (q^d-1) Y.$$
Furthermore, $\delta_q(Y) \subset Y^{\Gamma \#}$, and $Tr_q(Y^{\Gamma \#}) \subset Y^\#$.
\label{Lat}
\end{lem}
\proof
The inclusions are clear, since $n$ divides $q^d - 1$.  If $y \in Y$, and $y' \in Y^\Gamma$, then we find:
\begin{eqnarray*}
B(\delta_q y, y') & = & B(q \gamma y - y, y') \\
& = & q B(y, \gamma^{-1} y') - B(y, y') \\
& = & (q-1) B(y, y') \in n \ints, \mbox{ since } q-1 = mn.
\end{eqnarray*}
Hence $\delta_q(Y) \subset Y^{\Gamma \#}$.

Now, suppose that $w \in Y^{\Gamma \#}$, and $y' \in Y$.  Then, we find:
\begin{eqnarray*}
B(Tr_q(w), y') & = & \sum_{i=0}^{d-1} B(q^i \gamma^i w, y') \\
& \equiv & \sum_{i=0}^{d-1} B(\gamma^i w,  y') \Mod{n} \mbox{ since } q-1 = mn \\
& \equiv & B(w, Tr(y')) \in n \ints, \mbox{ since } Tr(y') \in Y^\Gamma.
\end{eqnarray*}
Hence $Tr_q(Y^{\Gamma \#}) \in Y^\#$.
\qed

Now, we fully describe $Z^\dag(T_F)$ with two results:
\begin{thm}
\label{Z1}
Suppose that $y_1 \in Y$, and $\bar y_2 \in \bar Y$.  Then if the element $t = y_1(\varpi) \bar y_2(\vartheta_L)$ is contained in $Z^\dag(T_F)$, then for every lift $y_2 \in Y$ of $\bar y_2$,
$$y_1, y_2 \in Y^\#, \mbox{ and } \delta_q y_2 \in (q^d-1) Y^{\Gamma \#}.$$
\end{thm}
\proof
For reference during this proof, we recall that:
$$nm = q-1, \mbox{ and } r = 1 + q + \cdots + q^{d-1}, \mbox{ and } nmr = q^d-1.$$
Suppose furthermore that $y_1', y_2' \in Y$, and let $\bar y_2' \in \bar Y$ be the reduction of $y_2'$.
Then, we find that $[Tr(y_1')](\varpi)$ and $[Tr_q(\bar y_2')](\vartheta_L)$ are elements of $T_F$.  It follows that:
\begin{eqnarray*}
C_L \left( y_1(\varpi) \bar y_2(\vartheta_L), [Tr(y_1')](\varpi) \right) & = & 1, \mbox{ and }\\
C_L \left( y_1(\varpi) \bar y_2(\vartheta_L), [Tr_q(\bar y_2')](\vartheta_L) \right) & = & 1.
\end{eqnarray*}

The explicit formula for the commutator $C_L$ yields:
\begin{eqnarray*}
1 & = & C_L \left(y_1(\varpi) y_2(\vartheta_L), [Tr_q(\bar y_2')](\vartheta_L) \right), \\
& = & \prod_{i=0}^{d-1} (\varpi, \vartheta_L)_{L,q^d-1}^{m q^i B(y_1, \gamma^i y_2')}, \\
& = & \zeta_L^{\sum_{i=0}^{d-1} m q^i B(y_1, \gamma^i y_2') }, \mbox{ since } (\varpi, \vartheta_L)_{L,q^d-1} = \zeta_L, \\
& = & \zeta_L^{\sum_{i=0}^{d-1} m q^i B(\gamma^{d-i} y_1, y_2')}, \mbox{ by the $\Gamma$-invariance of $B$}, \\
& = & \zeta_L^{\sum_{i=0}^{d-1} m q^i B(y_1, y_2')}, \mbox{ by the $\Gamma$-invariance of $y_1$}, \\
& = & \zeta_L^{ m r B(y_1, y_2') }, \mbox{ by summing a partial geometric series}, \\
& = & \zeta_F^{m B(y_1, y_2')}, \mbox{ since } \zeta_F = \zeta_L^r.
\end{eqnarray*}
Since $1 = \zeta_F^{m B(y_1, y_2')}$ for all $y_2' \in Y$, we find that:
$$y_1 \in Y^\#.$$

Carrying out a similar analysis, an explicit computation yields:
\begin{eqnarray*}
1 & = & C_L \left( y_1(\varpi) y_2(\vartheta_L), [Tr(y_1')](\varpi) \right), \\
& = & \prod_{i=0}^{d-1} (\varpi, \varpi)_{L, q^d-1}^{m B(y_1, \gamma^i y_1')} (\varpi, \vartheta_L)_{L, q^d-1}^{m B(y_2, \gamma^i y_1')}.
\end{eqnarray*}
Now, if $p$ is odd, we find that $q-1$ is even.  Since $(\varpi, \varpi)_{L,q^d-1} = \pm 1$, and $m B(y_1, \gamma^i y_1') \in mn \ints = (q-1) \ints \subset 2 \ints$ (since $y_1 \in Y^\#$), we find that:
$$(\varpi, \varpi)_{L,q^d-1}^{m B(y_1, \gamma^i y_1')} = 1.$$
On the other hand, if $p = 2$, $(\varpi, \varpi)_{L,q^d-1} = 1$, and once again the above equality holds.  Continuing our computations yields:
\begin{eqnarray*}
1 & = & \prod_{i=0}^{d-1} (\varpi, \varpi)_{L,q^d-1}^{m B(y_1, \gamma^i y_1')} (\varpi, \vartheta_L)_{L,q^d-1}^{m B(y_2, \gamma^i y_1')}, \\
& = & \prod_{i=0}^{d-1} \zeta_L^{m B(y_2, \gamma^i y_1')}, \mbox{ since } (\varpi, \varpi)_{L,q^d-1}^{m B(y_2, \gamma^i y_1')} = 1 \mbox{ and } (\varpi, \vartheta_L)_{L,q^d-1} = \zeta_L, \\
& = & \zeta_L^{\sum_{i=0}^{d-1} m q^i B(y_2, y_1') }, \mbox{ since } q \gamma(\bar y_2) = \bar y_2, \\
& = & \zeta_L^{m r B(y_2, y_1')}, \mbox{ by summing a partial geometric series,} \\
& = & \zeta_F^{m B(y_2, y_1')}, \mbox{ since } \zeta_F = \zeta_L^r.
\end{eqnarray*}
Hence, we find that $y_2 \in Y^\#$.

Finally, we prove that $\delta_q y_2 \in (q^d-1) Y^{\Gamma \#}$.  Note that $\delta_q y_2 \in (q^d - 1) Y$, since $\bar y_2 \in \bar Y^{\Gamma, q}$.  Thus, $\delta_q y_2 = (q^d-1) y_3$, for some $y_3 \in Y$.  It suffices to prove that $y_3 \in Y^{\Gamma \#}$.

Now, to prove that $y_3 \in Y^{\Gamma \#}$, suppose that $y' \in Y^\Gamma$.  It follows that:
\begin{eqnarray*}
1 & = & C_L\left( y_1(\varpi) \bar y_2(\vartheta_L), y'(\varpi) \right) \\
& = & (\varpi, \varpi)_{q^d-1}^{m B(y_1, y')} (\vartheta_L, \varpi)_{q^d-1}^{m B(y_2, y')} \\
& = & \zeta_L^{m B(y_2, y')}.
\end{eqnarray*}
Hence $B(y_2, y') \in nr \ints$.  It follows that:
\begin{eqnarray*}
B(y_3, y') & = & (q^d - 1)^{-1} B(\delta_q y_2, y') \\
& = & (q^d-1)^{-1} \left( B(q \gamma y_2, y') - B(y_2, y') \right) \\
& = & r^{-1} B(y_2, y') \in n \ints.
\end{eqnarray*}
Thus $y_3 \in Y^{\Gamma \#}$.
\qed

\begin{thm}
\label{Z2}
Suppose that $y_1, y_2 \in Y^\#$.  Also, suppose that $y_1 \in Y^\Gamma$, and $\bar y_2 \in \bar Y^{\Gamma, q}$.  Furthermore, suppose that $\delta_q y_2 \in (q^d-1)Y^{\Gamma \#}$.  Then $y_1(\varpi) \bar y_2(\vartheta_L) \in Z^\dag(T_F)$.
\end{thm}
\proof
Since $y_1 \in Y^\Gamma$, and $\bar y_2 \in \bar Y^{\Gamma, q}$, it follows that $y_1(\varpi) \bar y_2(\vartheta_L) \in T_F$.  Now, we may compute some commutators.

Suppose that $y_1' \in Y^\Gamma$, and $y_2' \in Y$, and $\bar y_2' \in \bar Y^{\Gamma, q}$.  Thus $y_1'(\varpi)$ and $\bar y_2'(\vartheta_L)$ are elements of $T_F$.  We begin by computing:
$$C_L(y_1(\varpi), y_1'(\varpi)) = (\varpi, \varpi)_{L, q^d-1}^{m B(y_1, y_1')}.$$
If $p$ is odd, then $mn = q-1$ is even, and thus $m B(y_1, y_1')$ is even.  Hence the commutator is trivial.  If $p$ is even, then $q^d-1$ is odd, and hence $(\varpi, \varpi)_{L, q^d-1} = 1$.  In either case, the commutator is trivial.

Now, consider the following commutator:
$$C_L(y_1(\varpi), \bar y_2'(\vartheta_L)) = \zeta_L^{m B (y_1, y_2')}.$$
We claim that $m B(y_1, y_2') \in (q^d-1) \ints$.  Indeed, we have:
$$B(y_1, y_2') = B(\gamma y_1, y_2') = B(y_1, \gamma^{-1} y_2') = B(y_1, q y_2' + (q^d-1) y_3'),$$
for some $y_3' \in Y$.  Since $y_1 \in Y^\#$, we have:
$$B(y_1, (q^d-1) y_3') \in n (q^d-1) \ints.$$
It follows that:
$$(q-1) B(y_1, y_2') \in n (q^d-1) \ints.$$
From this, we find:
$$B(y_1, y_2') \in nr \ints.$$
Hence $m B(y_1, y_2') \in mnr \ints = (q^d-1) \ints$.  This proves our claim, and we have proven that:
$$C_L(y_1(\varpi), \bar y_2'(\vartheta_L) = 1.$$

Next, consider the following commutator:
$$C_L(\bar y_2(\vartheta_L), y_1'(\varpi)) = \zeta_L^{-m B(y_2, y_1')}.$$
We claim now that $m B(y_2, y_1') \in (q^d-1) \ints$.  Indeed, we have:
$$B(y_2, y_1') = B(q \gamma y_2 + (q^d-1) y_3, y_1') = q B(y_2, y_1') + (q^d-1) B(y_3, y_1'),$$
for some $y_3 \in Y^{\Gamma \#}$.  In particular, $B(y_3, y_1') \in n \ints$, since $y_1' \in Y^\Gamma$.
It follows that:
$$(q-1) B(y_2, y_1') \in n (q^d-1) \ints.$$
From this we find that $B(y_2, y_1') \in nr \ints$, from which the claim follows.  We have proven that:
$$C_L(\bar y_2(\vartheta_L), y_1'(\varpi)) = 1.$$

Finally, note that $(\vartheta_L, \vartheta_L)_{L, q^d-1} = 1$.  Hence,
$$C_L(\bar y_2(\vartheta_L), \bar y_2'(\vartheta_L)) = 1.$$
We have proven that $y_1(\varpi)$, and $\bar y_2(\vartheta_L)$ commute with a set of generators for $T_F / T_F^1$.  Since $T_F^1 \in Z^\dag(T_F)$, this suffices to prove that:
$$y_1(\varpi) \bar y_2(\vartheta_L) \in Z^\dag(T_F).$$
\qed

The previous two theorems fully characterize the subgroup $Z^\dag(T_F)$.
\begin{cor}
\label{ZC}
Suppose that $y_1 \in Y$, $\bar y_2 \in \bar Y$, and $t^1 \in T_F^1$.  Then $t = y_1(\varpi) \bar y_2(\vartheta_L) t^1 \in Z^\dag(T_F)$ if and only if all of the following conditions hold:
\begin{itemize}
\item
$y_1 \in Y^{\# \Gamma}$.
\item
$y_2 \in Y^\#$, for any choice of representative $y_2$ of $\bar y_2$.
\item
$\delta_q y_2 \in (q^d-1) Y^{\Gamma \#}$, for any choice of representative $y_2$ of $\bar y_2$.
\end{itemize}
\end{cor}
\proof
This corollary follows directly from the previous two theorems.  One important observation is the following:  The latter two conditions do not depend upon the choice of representative $y_2 \in Y$ for a given $\bar y_2 \in \bar Y$.

Indeed, suppose that $y_2' = y_2 + (q^d-1) z$, for some $z \in Y$, so that $y_2$ and $y_2'$ are representatives for $\bar y_2$.  Since $Y^\# \subset n Y$, and $n$ divides $(q^d-1)$, we find that $y_2 \in Y^\#$ if and only if $y_2' \in Y^\#$.

Similarly, we find that $\delta_q y_2' = \delta_q y_2 + (q^d-1) \delta_q z$.  By Lemma \ref{Lat}, $\delta_q z \in Y^{\Gamma \#}$.  It follows that $\delta_q y_2 \in (q^d-1) Y^{\Gamma \#}$ if and only if $\delta_q y_2' \in (q^d-1) Y^{\Gamma \#}.$
\qed

The above corollary implies that $y_1(\varpi) \in Z^\dag(T_F)$, for a given $y_1 \in Y$, if and only if $y_1 \in Y^{\# \Gamma}$.  It also implies the following:
\begin{cor}
Suppose that $\bar y_2 \in \bar Y$.  Then $\bar y_2(\vartheta_L) \in Z^\dag(T_F)$ if and only if
\begin{equation}
\bar y_2 \in Im(Tr_q(\overline{Y^{\Gamma \#}}) \rightarrow \bar Y).
\label{EQa}
\end{equation}
\end{cor}
\proof
The previous corollary implies that $\bar y_2(\vartheta_L) \in Z^\dag(T_F)$ if and only if the following two conditions hold:
\begin{enumerate}
\item
$y_2 \in Y^\#$ for some (equivalently, every) representative $y_2$ of $\bar y_2$.
\item
$\delta_q y_2 \in (q^d-1) Y^{\Gamma \#}$ for some (equivalently, every) representative $y_2$ of $\bar y_2$.
\end{enumerate}

Given these conditions, and a representative $y_2$ of $\bar y_2$, there exists $w \in Y^{\Gamma \#}$ such that $\delta_q(y_2) = (q^d-1) w$.  Hence $\delta_q(y_2) = \delta_q Tr_q(w)$.  The injectivity of $\delta_q$ implies that $y_2 = Tr_q(w)$.  It follows that $\bar y_2$ is the image of $Tr_q(\bar w)$ in $\bar Y$.  Hence, the conditions (1) and (2) imply the one condition \ref{EQa} of this corollary.

Conversely, suppose that Equation \ref{EQa} is satisfied.  Then we may choose $w \in Y^{\Gamma \#}$, such that $\bar y_2$ equals the image of $Tr_q(\bar w)$ in $\bar Y$.  Thus $y_2 = Tr_q(w)$ is a representative for $\bar y_2$ in $\bar Y$.  Since $Tr_q(Y^{\Gamma \#}) \subset Y^\#$ by Lemma \ref{Lat}, the condition (1) is satisfied.  Since $\delta_q y_2 = Tr_q \delta_q w = (q^d-1) w$, condition (2) is satisfied as well.  Therefore, $\bar y_2(\vartheta_L) \in Z^\dag(T_F)$.
\qed

\subsection{The image of an isogeny}

For split metaplectic tori, we found a useful characterization of $Z^\dag(T_F)$ as the image of an isogeny on $F$-rational points.  The same isogeny makes sense for non-split tori, however there is a small but important difference between the image of the isogeny and $Z^\dag(T_F)$.  We view this difference as accounting for ``packets'' of representations of metaplectic tori, with the same parameter.

Consider the inclusion of $\ints[\Gamma]$-modules $\iota \colon Y^\# \hookrightarrow Y$.  Note that we use the fact that $Q$ is a $\Gamma$-invariant quadratic form, so that $Y^\#$ is a $\ints[\Gamma]$-submodule.  This inclusion corresponds to an isogeny of algebraic tori over $F$:
$$\alg{\iota} \colon \alg{T}^\# \rightarrow \alg{T}.$$
Our description of the $F$-rational and $L$-rational points for $\alg{T}$ is also valid, mutatis mutandis, for $\alg{T}^\#$.  When $y \in Y^\#$, and $u \in L^\times$, we simply write $(y \otimes u)$ for the corresponding element of $\alg{T}^\#(L) \ident Y^\# \otimes L^\times$.  We choose this notation, rather than $y(u)$, since we do not wish to confuse cocharacters of $\alg{T}$ with cocharacters of $\alg{T}^\#$.  Since $Y^\#$ is a $\ints[\Gamma]$-module, we find that:
\begin{prop}
The torus $\alg{T}^\#$ splits over an unramified extension of $F$.  Suppose that $y_1, y_2 \in Y^\#$.  Then $(y_1 \otimes \varpi) (y_2 \otimes \vartheta_L) \in T^\# = \alg{T}^\#(F)$ if and only if
$$y_1 \in Y^{\# \Gamma}, \mbox{ and } \bar y_2 \in (\overline{Y^\#})^{\Gamma, q}.$$
\end{prop}

The isogeny $\alg{\iota}$ has the following effect on $L$-rational points:
$$\iota(y \otimes u) = y(u), \mbox{ for all } y \in Y^\#, u \in L^\times, (y \otimes u) \in \alg{T}^\#(L).$$
Hence, we find that:
\begin{prop}
\label{II}
Suppose that $y_1 \in Y$, $\bar y_2 \in \bar Y$, and $t^1 \in T_L^1$.  Then $y_1(\varpi) \bar y_2(\vartheta) t^1$ is an element of the image of $\iota \colon \alg{T}^\#(F) \rightarrow \alg{T}(F)$ if and only if:
\begin{itemize}
\item
$y_1 \in Y^{\# \Gamma}.$
\item
There exists $y_2 \in Y^\#$ representing $\bar y_2$, such that $\delta_q y_2 \in (q^d-1) Y^\#$.
\item
$t^1 \in T_F^1$.
\end{itemize}
\end{prop}
\proof
Since $(n, p) = 1$, the image of $\iota$ contains $T_F^1$.  It suffices only to consider the images:
$$\iota((y_1 \otimes \varpi)(y_2 \otimes \vartheta_L)),$$
for all $y_1 \in Y^{\# \Gamma}$, and all $y_2 \in Y^\#$ such that $\bar y_2 \in (\overline{Y^\#})^{\Gamma, q}$.
\qed

Then, we find:
\begin{cor}
Suppose that $\bar y_2 \in \bar Y$.  Then $\bar y_2(\vartheta_L) \in \iota(\alg{T}^\#(F))$ if and only if
$$\bar y_2 \in Im(Tr_q(\overline{Y^\#}) \rightarrow \bar Y).$$
\end{cor}
\proof
If $\bar y_2 \in Im(Tr_q(\overline{Y^\#}) \rightarrow \bar Y)$, there exists an element $y_3 \in Y^\#$ such that $\bar y_2$ equals the image of $Tr_q(y_3)$ in $\bar Y$.  If $y_2 = Tr_q(y_3)$, then $y_2$ is a representative for $\bar y_2$ in $Y$.  Note that $y_2 \in Y^\#$, since $y_3 \in Y^\#$.  Furthermore, $\delta_q y_2 = \delta_q Tr_q(y_3) = (q^d-1) y_3 \in (q^d-1) Y^\#$.  Hence, $\bar y_2(\vartheta_L) \in \iota(\alg{T}^\#(F))$ by the previous proposition.

Conversely, suppose that $\bar y_2(\vartheta_L) \in \iota(\alg{T}^\#(F))$.  By the previous proposition, there exists a representative $y_2$ of $\bar y_2$ in $Y$, such that $y_2 \in Y^\#$, and $\delta_q(y_2) \in (q^d-1) Y^\#$.  It follows that $\delta_q(y_2) = \delta_q Tr_q(y_3)$, for some $y_3 \in Y^\#$.  Hence $y_2 = Tr_q(y_3)$.  Hence, $\bar y_2$ is contained in the image of $Tr_q(\overline{Y^\#})$ in $\bar Y$.
\qed

\subsection{The Packet Group}
From the previous two sections, we have described the groups $Z^\dag(T_F)$ and $\iota(\alg{T}^\#(F))$.  They are quite similar, with one exception.  Given $\bar y_2 \in \bar Y$, we have:
\begin{itemize}
\item
$\bar y_2(\vartheta_L) \in Z^\dag(T_F)$ if and only if
$$\bar y_2 \in Im(Tr_q(\overline{Y^{\Gamma \#}} ) \rightarrow \bar Y).$$
\item
$\bar y_2(\vartheta_L) \in \iota(\alg{T}^\#(F))$ if and only if
$$\bar y_2 \in Im(Tr_q(\overline{Y^\#}) \rightarrow \bar Y).$$
\end{itemize}

Define a finite group $P_{\theta_{\ff{l}}}^\dag$ by:
$$P_{\theta_{\ff{l}}}^\dag = { {Im(Tr_q(\overline{Y^{\Gamma \#}}) \rightarrow \bar Y) } \over {Im(Tr_q(\overline{Y^\#}) \rightarrow \bar Y)} }.$$
It follows from Proposition \ref{II}, and Corollary \ref{ZC}, that there is a natural short exact sequence:
\begin{equation}
1 \rightarrow \iota(\alg{T}^\#(F)) \rightarrow Z^\dag(T_F) \rightarrow P_{\theta_{\ff{l}}}^\dag \rightarrow 1.
\label{PES}
\end{equation}
However, this short exact sequence depends upon the choice of generator $\theta_{\ff{l}}$ of $\ff{l}^\times$.  We identify $P_{\theta_{\ff{l}}}^\dag$ here, in a way which makes the above short exact sequence independent of the choice of generator.

The $\ints[\Gamma]$-modules $Y^{\Gamma \#}$ and $Y^\#$ correspond to a pair of $\ff{f}$-tori $\alg{\bar T}^{\Gamma \#}$ and $\alg{\bar T}^\#$, which split over $\ff{l}$.  Moreover, the inclusions $Y^\# \subset Y^{\Gamma \#} \subset Y$ correspond to $\ff{f}$-isogenies of $\ff{f}$-tori:
$$\alg{\bar T}^\# \rightarrow \alg{\bar T}^{\Gamma \#} \rightarrow \alg{\bar T}.$$
The choice of generator $\theta_{\ff{l}}$ of $\ff{l}^\times$ corresponds to identifications:
$$\alg{\bar T}^{\Gamma \#}(\ff{l}) \ident \overline{Y^{\Gamma \#}}, \mbox{ and } \alg{\bar T}^\#(\ff{l}) \ident \overline{Y^\#}.$$
Furthermore, the trace map $Tr_q$ corresponds to the norm maps.  For example, there is a commutative diagram:
$$\xymatrix{ \alg{\bar T}^{\Gamma \#}(\ff{l}) \ar[r]^{\ident} \ar[d]^{N_{\ff{l}/\ff{k}}} & \overline{Y^{\Gamma \#}} \ar[d]^{Tr_q} \\
\alg{\bar T}^{\Gamma \#}(\ff{f}) \ar[r]^{\ident} & \overline{Y^{\Gamma \#}} }$$

Now, Lang's theorem \cite{Lang} implies that the norm map is surjective.  In other words, the commutative diagram above yields identifications:
$$\alg{\bar T}^{\Gamma \#}(\ff{f}) \ident Tr_q \overline{Y^{\Gamma \#}}, \mbox{ and } \alg{\bar T}^\#(\ff{f}) \ident Tr_q \overline{Y^\#}.$$
It follows that:
\begin{prop}
Define a finite group $P^\dag$ by:
$$P^\dag = { {Im(\alg{\bar T}^{\Gamma \#}(\ff{f}) \rightarrow \alg{\bar T}(\ff{f}))} \over {Im(\alg{\bar T}^{\#}(\ff{f}) \rightarrow \alg{\bar T}(\ff{f}))} }.$$
Then $P_{\theta_{\ff{l}}}^\dag \isom P^\dag$.
\end{prop}
The construction of the group $P^\dag$ does not depend upon the choice of generator $\theta_{\ff{l}}$. This yields a canonical short exact sequence, which does not depend upon the choice of uniformizing element $\varpi$ nor upon the choice of generator $\theta_{\ff{l}}$:
$$1 \rightarrow \iota(\alg{T}^\#(F)) \rightarrow Z^\dag(T_F) \rightarrow P^\dag \rightarrow 1.$$

The main theorem of Langlands in \cite{Lan}, parameterizing smooth characters of tori over local fields, determines isomorphisms:
$$\Sp{X}(T_F) \isom H_c^1(\Weil_{L/F}, \Sp{\hat T}), \mbox{ and } \Sp{X}(T_F^\#) \isom H_c^1(\Weil_{L/F}, \Sp{\hat T}^\#).$$
As before, the characters of the image of an isogeny can be parameterized cohomologically:
\begin{prop}
The Langlands parameterization yields a {\em finite-to-one} parameterization of the smooth characters of $Z^\dag(T_F)$:
$$1 \rightarrow \Sp{X}(P^\dag) \rightarrow \Sp{X}(Z^\dag(T_F)) \rightarrow { {H_c^1(\Weil_{L/F}, \Sp{\hat T})} \over {H_c^1(\Weil_{L/F}, \Sp{\hat T} \rightarrow \Sp{\hat T}^\#)} } \rightarrow 1.$$
\end{prop}
\begin{rem}
In order to view $H_c^1(\Weil_{L/F}, \Sp{\hat T} \rightarrow \Sp{\hat T}^\#)$ as a subgroup of $H_c^1(\Weil_{L/F}, \Sp{\hat T})$ as above, we must know that the map $H_c^0(\Weil_{L/F}, \Sp{\hat T}) \rightarrow H_c^0(\Weil_{L/F}, \Sp{\hat T^\#})$ is surjective.  This follows from the identifications:
$$H_c^0(\Weil_{L/F}, \Sp{\hat T}) \ident \Hom_\ints(Y^\Gamma, \complex^\times),  \mbox{ and } H_c^0(\Weil_{L/F}, \Sp{\hat T}^\#) \ident \Hom_\ints(Y^{\# \Gamma}, \complex^\times),$$
and the fact that $Y^{\# \Gamma}$ has finite index in $Y^\Gamma$.
\end{rem}

This leads directly, via a Stone-von-Neumann theorem, to a main theorem for tame covers of unramified tori:
\begin{thm}
\label{TT}
Suppose that we are given a {\em tame} metaplectic cover of an {\em unramified} torus:
$$1 \rightarrow \mu_n \rightarrow \tilde T \rightarrow T \rightarrow 1.$$
Then, with the sublattices $Y^\# \subset Y^{\Gamma \#} \subset Y$ defined as before, and the resulting isogenies $\alg{T}^\# \rightarrow \alg{T}^{\Gamma \#} \rightarrow \alg{T}$ of unramified tori, we find:
\begin{quote}
There is a finite-to-one surjective map, intertwining the natural action of $H^1(\Weil_{L/F}, \Sp{\hat T})$:
$$\Phi \colon \Sp{Irr}_\epsilon(\tilde T) \rightarrow  { {H_c^1(\Weil_{L/F}, \Sp{\hat T})} \over {H_c^1(\Weil_{L/F}, \Sp{\hat T} \rightarrow \Sp{\hat T}^\#)} }.$$
The fibres of this map are torsors for the finite group $\Sp{X}(P^\dag)$, where:
$$P^\dag = { {Im(\alg{\bar T}^{\Gamma \#}(\ff{f}) \rightarrow \alg{\bar T}(\ff{f}))} \over {Im(\alg{\bar T}^{\#}(\ff{f}) \rightarrow \alg{\bar T}(\ff{f}))} }.$$
\end{quote}
\end{thm}

\begin{rem}
We do not know if a parameterization, such as that above, holds for general metaplectic tori over local fields.  Namely, we have not been able to describe the center of such metaplectic tori, when $\alg{T}$ is ramified, or when $\alg{T}$ is an unramified torus, but the cover is not tame.  We hope that such a parameterization is possible, though the packets might be substantially different.
\end{rem}

\begin{rem}
In the process of proving the previous theorem, we chose a uniformizing element $\varpi \in F^\times$ and a root of unity $\theta_L$.  However, this choice does not have any effect on the parameterization given above.  The sublattices $Y^\#$ and $Y^{\Gamma \#}$ clearly do not depend upon such a choice.  Moreover, the action of $\Sp{X}(P^\dag)$ on the fibres of $\Phi$ does not depend on such a choice.
\end{rem}

\section{Pseudo-Spherical and Pseudo-Trivial Representations}
\label{PS}
We maintain all of the conventions of the previous section.  In particular, we have a {\em tame} metaplectic cover of an {\em unramified} torus:
$$1 \rightarrow \mu_n \rightarrow \tilde T_F \rightarrow T_F \rightarrow 1.$$
We have shown that the irreducible genuine representations of $\tilde T_F$ can be parameterized by the points of a homogeneous space on which $H^1(\Weil_{L/F}, \Sp{\hat T})$ acts transitively.  However, such a parameterization is {\em not} unique; one must choose a ``base point'' in the space of irreducible genuine representations of $\tilde T_F$, in order to choose a specific morphism of homogeneous spaces:
$$\Phi \colon \Sp{Irr}_\epsilon(\tilde T_F) \rightarrow { {H^1(\Weil_{L/F}, \Sp{\hat T})} \over {H^1(\Weil_{L/F}, \Sp{\hat T} \rightarrow \Sp{\hat T}^\#)} }.$$

In this section, we discuss the data which determines such base points.  Such choices arise frequently in treatments of metaplectic groups, often as choices of square roots of $-1$ in $\complex$.

\subsection{The Residual Extension}

Recall that the unramified torus $\alg{T}$ has a smooth model $\scheme{T}$ over $\OO_F$, and $T_F^\circ = \scheme{T}(\OO_F)$.  In this case, $T_F^\circ$ is the maximal compact subgroup of $T$, and we let $\tilde T_F^\circ$ be its preimage in $\tilde T_F$.  Also, $\alg{\bar T}$ denotes the special fibre of $\scheme{T}$, which is a torus over $\ff{f}$.  Recall that $\alg{T}'$ is a central extension of $\alg{T}$ by $\alg{K}_2$.  Pushing forward via the tame symbol led to the tame central extension:
$$1 \rightarrow \ff{f}^\times \rightarrow T_F^t \rightarrow T_F \rightarrow 1.$$
We write $T_F^{t \circ}$ for the preimage of $T_F^\circ$ in $T_F^t$.

In Section 12.11 of \cite{D-B}, Deligne and Brylinski construct an extension $\alg{\bar T}'$ of $\alg{\bar T}$ by $\alg{G}_m$ (in the category of groups over $\ff{f}$).  We call $\alg{\bar T}'$ the {\em residual extension} associated to $\alg{T}'$.  The residual extension fits into the following commutative diagram:
$$\xymatrix{ 1 \ar[r] & \ff{f}^\times \ar[r] \ar[d] & T_F^{t \circ} \ar[r] \ar[d] & T_F^\circ \ar[r] \ar[d] & 1 \\
1 \ar[r] & \ff{f}^\times \ar[r] & \bar T_\ff{f}' \ar[r] & \bar T_\ff{f} \ar[r] & 1.}$$
Here, the map from $\ff{f}^\times$ to itself is the identity, the map from $T_F^\circ$ to $\bar T_\ff{f}$ is the reduction map, and the diagram identifies the top row with the pullback of the bottom row via reduction.

As an extension of $\alg{\bar T}$ by $\alg{G}_m$ over $\ff{f}$, the group $\alg{\bar T}'$ is an algebraic torus over $\ff{f}$.  Note that the category of extensions of $\alg{\bar T}$ by $\alg{G}_m$, in the category of groups over $\ff{f}$, is equivalent to the category of extensions of $Y$ by $\ints$, in the category of $\ints[\Gamma]$-modules (where $\ints$ is given the trivial module structure).  In this way, the construction of Section 12.11 of \cite{D-B} associates an extension of $Y$ by $\ints$, to any extension of an unramified torus $\alg{T}$ by $\alg{K}_2$.
\begin{rem}
Recall that $\tilde Y$ is a $\Gamma$-equivariant extension of $Y$ by $L^\times$, constructed as a functorial invariant of the extension $\alg{T}'$ of $\alg{T}$ by $\alg{K}_2$.  Let $Y'$ be the extension of $Y$ by $\ints$, obtained by pushing forward $\tilde Y$ via the valuation map $L^\times \rightarrow \ints$:
$$0 \rightarrow \ints \rightarrow Y' \rightarrow Y \rightarrow 0.$$
We do not know whether this extension is naturally isomorphic to the exact sequence of cocharacter groups of the residual extension of tori described above
\end{rem}

\begin{defn}
Let $Spl(\alg{\bar T}')$ denote the set of splittings, in the category of algebraic groups over $\ff{f}$, of the short exact sequence:
$$1 \rightarrow \alg{G}_m \rightarrow \alg{\bar T}' \rightarrow \alg{\bar T} \rightarrow 1.$$
We say that the extension $\alg{T}'$ of the unramified torus $\alg{T}$ is a {\em residually split extension}, if $Spl(\alg{\bar T}')$ is non-empty.
\end{defn}
In particular, if $\alg{T}$ is a split torus, then $\alg{T}'$ is residually split.

\begin{prop}
If $Spl(\alg{\bar T}')$ is non-empty, then $Spl(\alg{\bar T}')$ is a torsor for the abelian group $X^\Gamma$.
\end{prop}
\proof
Any two algebraic splittings are related by an element of $\Hom_\ff{f}(\alg{\bar T}, \alg{G}_m)$.  This group may be identified with the $\Gamma$-fixed characters of $\alg{T}$.
\qed

\subsection{Pseudo-spherical representations}

Suppose now that $\alg{T}'$ is a residually split extension of $\alg{T}$ by $\alg{K}_2$.  Fix a splitting $s \in Spl(\alg{\bar T}')$.  The splitting lifts to a splitting $\sigma \colon T_F^\circ \rightarrow T_F^{t \circ}$.  Pushing forward via the $m^{th}$ power map, we may also view $\sigma$ as a splitting $T_F^\circ \rightarrow \tilde T_F^\circ$.  From such a splitting $s$, we let $\theta_s^\circ \colon \tilde T_F^\circ \rightarrow \complex^\times$ denote the character obtained by projecting onto $\mu_n$ (via the splitting $\sigma$), and then applying the injective homomorphism $\epsilon \colon \mu_n \rightarrow \complex^\times$.

Let $Z_{\tilde T_F}(\tilde T_F^\circ)$ be the centralizer of $\tilde T_F^\circ$ in $\tilde T_F$.  Then, we find:
\begin{prop}
The group $Z_{\tilde T_F}(\tilde T_F^\circ)$ is the preimage of a subgroup $Z_{T_F}^\dag(T_F^\circ) \subset T_F$.  Considering the valuation map:
$$val \colon T_F \rightarrow Y^\Gamma,$$
whose kernel is $T_F^\circ$, $Z_{T_F}^\dag(T_F^\circ)$ is equal to the preimage of $Y^{\# \Gamma }$.
\end{prop}
\proof
Since $Z_{T_F}^\dag(T_F^\circ) \supset T_F^\circ$, it suffices to identify the set of $y \in Y^\Gamma$ such that:
$$C_L(y(\varpi), \bar y'(\vartheta_L)) = 1 \mbox{ for all } \bar y' \in \bar Y^{\Gamma, q}.$$
In fact, the set of such $y$ has been identified in the proofs of Theorems \ref{Z1} and \ref{Z2}.  The above condition is satisfied if and only if $y \in Y^{\# \Gamma}$.
\qed

\begin{cor}
The group $Z_{\tilde T_F}(\tilde T_F^\circ)$ is abelian.
\end{cor}
\proof
As $Z_{\tilde T_F}(\tilde T_F^\circ)$ is the centralizer of the abelian group $\tilde T_F^\circ$, it suffices to prove that $C(y(\varpi), y'(\varpi)) = 1$ for all $y,y' \in Y^{\# \Gamma}$.  This is proven in the beginning of the proof of Theorem \ref{Z2}.
\qed

Directly following Section 4 of Savin \cite{Sav}, we find
\begin{prop}
There is a natural bijection between the following two sets:
\begin{itemize}
\item
The set $\Sp{Irr}_{s,\epsilon}^{sph}(\tilde T_F)$ of {\em pseudo-spherical} irreducible representations of $\tilde T_F$ (for the splitting $s$).  These are the genuine irreducible representations of $\tilde T_F$, whose restriction to $\tilde T_F^\circ$ via the splitting $s$ contains a nontrivial $\theta_s^\circ$-isotypic component.
\item
The set of {\em extensions} of $\theta_s^\circ$ to the group $Z_{\tilde T_F}(\tilde T_F^\circ)$.
\end{itemize}
Namely, if $(\pi, V)$ is a pseudo-spherical irreducible representation, its $\theta_s^\circ$-isotypic subrepresentation is an extension of $\theta_s^\circ$ to the group $Z_{\tilde T_F}(\tilde T_F^\circ)$.  Conversely, given such an extension $\theta_s^1$ of $\theta_s^\circ$ to a character of $Z_{\tilde T_F}(\tilde T_F^\circ)$, the induced representation $Ind_{Z_{\tilde T_F}(\tilde T_F^\circ)}^{\tilde T_F} \theta_s^1$ is a pseudo-spherical irreducible representation.
\end{prop}

One may rephrase the above bijection slightly; the splitting $s$ yields an injective homomorphism from $T_F^\circ$ onto a normal subgroup of $Z_{\tilde T_F}(\tilde T_F^\circ)$.  This fits into a commutative diagram with exact rows and columns:
$$\xymatrix{& & 1 \ar[d] & 1 \ar[d] & \\
1 \ar[r] & 1 \ar[r] \ar[d] & T_F^\circ \ar[r] \ar[d]^{s} & T_F^\circ \ar[r] \ar[d] & 1 \\
1 \ar[r] & \mu_n \ar[r] \ar[d] & Z_{\tilde T_F}(\tilde T_F^\circ) \ar[r] \ar[d] & Z_{T_F}^\dag(T_F^\circ) \ar[r] \ar[d] & 1 \\
1 \ar[r] & \mu_n \ar[r] \ar[d] & \tilde Y^{\# \Gamma} \ar[r] \ar[d]  & Y^{\# \Gamma} \ar[r] \ar[d] & 1 \\
& 1 & 1 & 1 &
}$$
Hence, the splitting $s$ determines an extension $\tilde Y^{\# \Gamma}$ of $Y^{\# \Gamma}$ by $\mu_n$.  A standard diagram chase now yields:

\begin{prop}
There is a natural bijection:
$$\Sp{Irr}_{s,\epsilon}^{sph}(\tilde T_F) \leftrightarrow \Sp{X}_\epsilon(\tilde Y^{\# \Gamma}).$$
\end{prop}

\begin{cor}
The space $\Sp{Irr}_{s,\epsilon}^{sph}(\tilde T_F)$ is naturally a torsor for the complex algebraic torus $\Sp{X}(Y^{\# \Gamma})$.
\end{cor}

\begin{rem}
One may view $\Sp{X}_\epsilon(\tilde Y^{\# \Gamma})$ as the set of irreducible representations of a ``quantum torus''.  Indeed, consider the ring:
$$\complex_\epsilon[\tilde Y^{\# \Gamma}] = { {\complex[\tilde Y^{\# \Gamma}]} \over {\langle \zeta - \epsilon(\zeta) \rangle_{\zeta \in \mu_n}} }.$$
The ring $\complex_\epsilon[\tilde Y^{\# \Gamma}]$ can be viewed as (the coordinate ring of) a quantum complex torus, which we call $\Sp{\hat T}_\epsilon^{\# \Gamma}$.  $\Sp{\hat T}_\epsilon^{\# \Gamma}$ is the quantization of a complex torus, at a root of unity.  Quasi-coherent sheaves on this quantum torus (i.e. modules over its coordinate ring) correspond naturally to pseudo-spherical representations of $\tilde T_F$.
\end{rem}

\subsection{Pseudo-Trivial Representations}

In many practical situations, the extension $\tilde Y^{\# \Gamma}$ of $Y^{\# \Gamma}$ by $\mu_n$ splits over a quite large submodule of $Y^{\# \Gamma}$.  For example, in many cases, the extension splits over $Y^{\# \Gamma} \cap 2 Y$.

Suppose that $V \subset Y^{\# \Gamma}$ is a finite index subgroup, endowed with a splitting $v$ of the resulting exact sequence:
$$\xymatrix{
1 \ar[r] & \mu_n \ar[r] & \tilde V \ar[r] & V \ar[r] \ar @/_/ [l]_{v}& 1.}$$
Let $U = Y^{\# \Gamma} / V$ denote the quotient.  The splitting $v$ yields an extension of finite abelian groups:
$$1 \rightarrow \mu_n \rightarrow \tilde U \rightarrow U \rightarrow 1.$$

Pulling back yields natural inclusions:
$$\Sp{X}_\epsilon(\tilde U) \hookrightarrow \Sp{X}_\epsilon(\tilde Y^{\# \Gamma}) \ident \Sp{Irr}_{s,\epsilon}^{sph}(\tilde T_F).$$

Therefore, within the set of pseudo-spherical representations of $\tilde T_F$, we find a {\em finite} set of ``pseudo-trivial'' representations (relative to the choice of splitting subgroup $(V,v)$ of $Y^{\# \Gamma}$:
\begin{defn}
The genuine pseudo-trivial representations of $\tilde T_F$, are those irreducible pseudo-spherical genuine representations, that are in the image of $\Sp{X}_\epsilon(\tilde U)$.  This definition depends upon the following choices:
\begin{itemize}
\item
The splitting $s$ (to determine the pseudo-spherical representations).
\item
The splitting subgroup $(V,v)$ (to determine the pseudo-trivial representations).
\end{itemize}
\end{defn}

\begin{rem}
Most often, one chooses a pseudo-trivial ``base point'' in the space $\Sp{Irr}_\epsilon(\tilde T_F)$.  Very often (cf. the examples of \cite{Sav}) $\tilde U$ is a finite abelian group of exponent $4$.  It follows that pseudo-trivial representations may often be given by specifying certain characters of an abelian group of exponent $4$.  This explains the frequent need to choose fourth roots of unity, in the literature on metaplectic groups.
\end{rem}

\section{Tori over $\reals$}

In this section, the following will be fixed:
\begin{itemize}
\item
$\alg{T}$ will be a torus over $\reals$, and $\Gamma = Gal(\complex/\reals) = \{ 1, \gamma \}$.  $X$ and $Y$ will be the resulting character and cocharacter groups, viewed as $\ints[\Gamma]$-modules.
\item
$\alg{T}'$ will be an extension of $\alg{T}$ by $\alg{K}_2$ in $\Cat{Gp}_\reals$.
\item
$(Q,\tilde Y)$ will be the Deligne-Brylisnki invariants of $\alg{T}'$.  $B$ will be the symmetric bilinear form associated to $Q$.
\item
We fix $n = 2$, so that $\reals$ has enough $n^{th}$ roots of unity.
\item
If $W$ is a subgroup of $Y$, then we write $W^\Gamma$ for the subgroup of $\Gamma$-fixed elements of $W$.  We also define:
$$W^\# = \{ y \in Y \mbox{ such that } B(y, w) \in 2 \ints \mbox{ for all } w \in W \}.$$
\item
$\epsilon \colon \alg{\mu}_2(\reals) \rightarrow \complex^\times$ will be the unique injective character.
\item
We view $T = \alg{T}(\reals)$ as a real Lie group.  We write $T^\circ$ for the connected component of the identity element, and $\pi_0 T$ for the component group of $T$.
\item
The extension $\alg{T}'$, and the quadratic Hilbert symbol,  yields an extension of Lie groups:
$$1 \rightarrow \mu_2 \rightarrow \tilde T \rightarrow T \rightarrow 1.$$
\end{itemize}

We are interested in parameterizing the irreducible genuine representations of $\tilde T$, and the set of such representations is called $\Sp{Irr}_\epsilon(\tilde T)$, as before.

\subsection{Structure of metaplectic tori over $\reals$}
There is a short exact sequence of Lie groups:
$$1 \rightarrow T^\circ \rightarrow T \rightarrow \pi_0 T \rightarrow 1.$$
Let $\alg{T}^\Gamma$ be the split real torus with cocharacter group $Y^\Gamma$.  Let $T^\Gamma$ denote the real points of $\alg{T}^\Gamma$.  Then, we find that $\pi_0 T^\Gamma$ is canonically isomorphic to $\overline{Y^\Gamma} = Y^\Gamma \otimes_\ints \mu_2 \ident Y^\Gamma / 2 Y^\Gamma$.  Moreover, the inclusion of $\reals$-tori from $\alg{T}^\Gamma$ into $\alg{T}$ induces a surjective map of component groups:
$$\pi_0 T^\Gamma \twoheadrightarrow \pi_0 T.$$
Therefore, every element $t$ of $T$ has a (often non-unique) decomposition $t = t^\circ \bar y(-1)$, for some $t^\circ \in T^\circ$ and $\bar y \in \overline{Y^\Gamma}$.  In other words, there is a natural surjective homomorphism:
$$\overline{Y^\Gamma} \rightarrow \pi_0 T.$$

Now, we consider the metaplectic cover of $T$:
$$1 \rightarrow \mu_2 \rightarrow \tilde T \rightarrow T \rightarrow 1.$$
The commutator $C \colon T \times T \rightarrow \mu_2$ is bi-multiplicative and continuous.  It follows that the commutator is trivial when either of its inputs is in $T^\circ$.  Hence we find:
\begin{prop}
$T^\circ$ is a subgroup of $Z^\dag(T)$.
\end{prop}

\subsection{Description of the center}

It follows from the previous proposition that, to describe $Z^\dag(T)$, it suffices to describe its image in $T / T^\circ$.  Hence, it suffices to determine for which $\bar y \in \overline{Y^\Gamma}$, it holds that $\bar y(-1) \in Z^\dag(T)$.  We must be able to compute the commutator $C(\bar y(-1), \bar y'(-1))$ for arbitrary $\bar y, \bar y' \in \overline{Y^\Gamma}$.

Here, we note that such elements $\bar y(-1)$ and $\bar y'(-1)$ are contained in the real points of the maximal $\reals$-split subtorus $\alg{T}^\Gamma \hookrightarrow \alg{T}$.  Restricting the central extension of $\alg{T}$ by $\alg{K}_2$, to the split subtorus $\alg{T}^\Gamma$, the formula of Corollary 3.14 of \cite{D-B} is valid for computing commutators.  We find that:
\begin{prop}
If $y, y' \in Y^\Gamma$, then $C(\bar y(-1), \bar y'(-1)) = (-1)^{B(y, y')}$.
\end{prop}
\proof
This follows directly from Corollary 3.14 of \cite{D-B}, and the Hilbert symbol over $\reals$:  $(-1, -1)_{\reals, 2} = -1$.
\qed

\begin{prop}
Given $\bar y \in \overline{Y^\Gamma}$, $\bar y(-1) \in Z^\dag(T)$ if and only if every representative $y$ of $\bar y$ in $Y$ satisfies:
$$y \in Y^{\Gamma \# \Gamma}.$$
\end{prop}
\proof
Suppose $\bar y, \bar y' \in \overline{Y^\Gamma}$.  Let $y$ be a representative of $\bar y$ in $Y$.  The commutator has been computed:
$$C(\bar y(-1), \bar y'(-1)) = (-1)^{B(\bar y, \bar y')}.$$
Thus, we find that $C(\bar y(-1), \bar y'(-1)) = 1$ for all $\bar y' \in \overline{Y^\Gamma}$, if and only if $B(y, y') \in 2 \ints$ for all represenatives $y$ of all $\bar y' \in \overline{Y^\Gamma}$.  This occurs if and only if $B(y, y') \in 2 \ints$ for all $y' \in Y^\Gamma$, i.e., $y \in Y^{\Gamma \#}$.

Thus, we find that, given $\bar y \in \overline{Y^\Gamma}$, $\bar y(-1) \in Z^\dag(T)$ if and only if $y \in Y^{\Gamma \#} \cap Y^\Gamma = Y^{\Gamma \# \Gamma}$.
\qed

\begin{cor}
Let $\alg{T}^{\Gamma \#}$ be the real torus with cocharacter group $Y^{\Gamma \#}$.  Let $T^{\Gamma \#} = \alg{T}^{\Gamma \#}(\reals)$.  Then, the quotient $Z^\dag(T) / T^\circ$ is isomorphic to $Im(\pi_0 T^{\Gamma \#} \rightarrow \pi_0 T)$.
\end{cor}
\proof
There is a commutative diagram of finite abelian groups:
$$\xymatrix{
Y^{\Gamma \# \Gamma} \otimes_\ints \mu_2 \ar[r]^{\eta} \ar[d]^{p^{\Gamma \#}} & Y^\Gamma \otimes_\ints \mu_2 \ar[d]^{p} \\
\pi_0 T^{\Gamma \#} \ar[r]^{\rho} & \pi_0 T. }$$
The previous proposition demonstrates that $Z^\dag(T) / T^\circ$ can be identified with the image of $p \circ \eta$.  The commutativity of the above diagram, together with the surjectivity of $p^{\Gamma \#}$ implies that this image is the same as the image of $\rho$.
\qed

\subsection{The image of an isogeny}
As in the nonarchimedean case, the inclusion $Y^\# \hookrightarrow Y$ of $\ints[\Gamma]$-modules corresponds to an isogeny of tori over $\reals$:
$$\alg{\iota} \colon \alg{T}^\# \rightarrow \alg{T}.$$
We are interested in the resulting continuous homomorphism of real Lie groups:
$$\iota \colon T^\# \rightarrow T.$$
Since $\alg{\iota}$ is an isogeny, we find that $\iota(T^\#) \supset T^\circ$.  Thus, in order to fully describe $\iota(T^\#)$ it suffices to determine for which $\bar y \in \overline{Y^\Gamma}$, $\bar y(-1) \in \iota(T^\#)$.
\begin{prop}
Suppose that $\bar y \in \overline{Y^\Gamma}$.  Then $\bar y(-1) \in \iota(T^\#)$ if and only if $\bar y \in Im( \overline{Y^{\# \Gamma}} \rightarrow \overline{Y^\Gamma})$.
\end{prop}
\proof
We find that $\bar y(-1) \in \iota(T^\#)$, if and only if there exists $y \in Y^{\# \Gamma}$ which represents $\bar y$.  The proposition follows.
\qed

\begin{cor}
The quotient $\iota(T^\#) / T^\circ$ can be identified with the image $Im(\pi_0 T^\# \rightarrow \pi_0 T)$.
\end{cor}
Comparing the image of the isogeny $\iota$, to the group $Z^\dag(T)$, yields a short exact sequence:
$$1 \rightarrow \iota(T^\#) \rightarrow Z^\dag(T) \rightarrow P^\dag \rightarrow 1,$$
where we may identify the finite group:
$$P^\dag \ident { {Im(\pi_0 T^{\Gamma \#} \rightarrow \pi_0 T)} \over {Im(\pi_0 T^\# \rightarrow \pi_0 T)} }.$$

\subsection{Parameterization}
As for the case of nonarchimedean fields, we choose to parameterize the genuine irreducible representations of $\tilde T$, through a finite-to-one map and a description of the fibres.  Over $\reals$, the previous two sections imply that the space $\Sp{Irr}_\epsilon(\tilde T)$ can be identified (via Theorem \ref{SvN}) with the complex variety of genuine characters $\Sp{X}_\epsilon(Z(\tilde T))$.  This is a torsor for the complex algebraic group of characters $\Sp{X}(Z^\dag(T))$.  There is a short exact sequence:
$$1 \rightarrow \Sp{X}(P^\dag) \rightarrow \Sp{X}(Z^\dag(T)) \rightarrow { {H_c^1(\Weil_\reals, \Sp{\hat T})} \over {H_c^1(\Weil_\reals, \Sp{\hat T} \rightarrow \Sp{\hat T}^\#)} } \rightarrow 1.$$
Hence, we find:
\begin{thm}
\label{RT}
Suppose that we are given a metaplectic cover of a real torus:
$$1 \rightarrow \mu_n \rightarrow \tilde T \rightarrow T \rightarrow 1.$$
Then, with the sublattices $Y^\# \subset Y^{\Gamma \#} \subset Y$ defined as before, and the resulting isogenies $\alg{T}^\# \rightarrow \alg{T}^{\Gamma \#} \rightarrow \alg{T}$, we find:
\begin{quote}
There is a finite-to-one surjective map, intertwining the natural action of $H^1(\Weil_\reals, \Sp{\hat T})$:
$$\Phi \colon \Sp{Irr}_\epsilon(\tilde T) \rightarrow  { {H_c^1(\Weil_\reals, \Sp{\hat T})} \over {H_c^1(\Weil_\reals, \Sp{\hat T} \rightarrow \Sp{\hat T}^\#)} }.$$
The fibres of this map are torsors for the finite group $\Sp{X}(P^\dag)$, where:
$$P^\dag = { {Im(\pi_0 T^{\Gamma \#} \rightarrow \pi_0 T)} \over {Im(\pi_0 T^\# \rightarrow \pi_0 T)} }.$$
\end{quote}
\end{thm}

Note that this theorem is quite similar to the parameterization of $\Sp{Irr}_\epsilon(\tilde T)$ for tame covers of unramified tori over nonarchimedean local fields.  The only difference is that points of residual tori are replaced by component groups.

%----------------------------------------------------------------
\bibliographystyle{amsplain}
\bibliography{MetaplecticTori}
\end{document}